\documentclass[12pt]{amsart}
\usepackage{mathrsfs}

\usepackage{amssymb,amsfonts,amsthm,amsmath}
\usepackage{epsfig}
\usepackage[utf8]{inputenc}
\usepackage{mathrsfs}
\usepackage{lscape}
\usepackage{amssymb,amsmath,graphicx,color,textcomp, amsthm,bbm,bbold, enumerate,booktabs}
\usepackage{chngcntr}
\usepackage{apptools}
\AtAppendix{\counterwithin{theorem}{section}}
\definecolor{Red}{cmyk}{0,1,1,0}

\definecolor{verde}{cmyk}{1,0,1,0}

\definecolor{loka}{cmyk}{.5,0,1,.5}
\definecolor{azul}{cmyk}{1,1,0,0}

\numberwithin{equation}{section}

\newcommand{\be}{\begin{equation}}
\newcommand{\ee}{\end{equation}}

\newtheorem{definition}{Definition}

\newtheorem{corolario}{Corollary}

\newtheorem{lemma}{Lemma}
\newtheorem{teorema}{Theorem}

\begin{document}
\title{A note on the mild solutions of Hilfer impulsive fractional differential equations}
\author{J. Vanterler da C. Sousa$^1$}
\address{$^1$ Department of Applied Mathematics, Institute of Mathematics,
 Statistics and Scientific Computation, University of Campinas --
UNICAMP, rua S\'ergio Buarque de Holanda 651,
13083--859, Campinas SP, Brazil\newline
$^2$ Coordination of Civil Engineering, Technological Federal University of Paran\'{a}, 85053-525, Guarapuava/PR, Brazil
e-mail: {\itshape \texttt{vanterlermatematico@hotmail.com, oliveiradaniela@utfpr.edu.br, capelas@ime.unicamp.br }}}
\author{D. S. Oliveira$^2$}
\author{E. Capelas de Oliveira$^1$}

\begin{abstract} In this paper, we present a new class of inequality of the Gronwall type and discuss some particular cases. In this sense, we investigate the uniqueness and $\delta$-Ulam-Hyers-Rassias stability of the mild solution for the fractional differential equation with non-instantaneous impulses in a P$\delta$-normed Banach space.

\newline 
\end{abstract}
\maketitle

\section{Introduction}\label{intro}

Over time, the study of differential equations has been object of study by numerous researchers \cite{braun,farlow}, from the practical context, in which they investigated through mathematical models that describe natural phenomena and in the theoretical sense, such as the investigation of the existence, uniquenesses, Ulam-Hyers stability, among others, of classic and mild solutions of differential equations such as functional, impulsive, evolution, sometimes with non-instantaneous or instantaneous impulses, among others \cite{chen,farlow,hernandez,princ,shi1}. The range of options that the field of differential equations enabled the researchers is in an exponential growth that countless works over time were published \cite{andras2013,bala,chen,jungu,radha,rassias}.

In 2013, Hernandez and O'Regan \cite{hernandez} investigated the existence of a new class of impulsive differential equations with non-instantaneous impulses, given by the problem
\begin{equation}
\left\{ 
\begin{array}{cll}
u^{\prime }\left( t\right)  & = & Au\left( t\right) +f\left( t,u\left(
t\right) \right) ,\text{ }t\in \left( s_{i},t_{i+1}\right] ,\text{ }i=0,...,N
\\ 
u\left( t\right)  & = & \gamma _{k}\left( t,u\left( t\right) \right) ,\text{ 
}t\in \left( t_{i},s_{i}\right] ,\text{ }i=1,...,N \\ 
u\left( 0\right)  & = & x_{0}%
\end{array}%
\right. 
\end{equation}
where $A: D(A)\subset X \rightarrow X$ is the generator of $C_{0}$-semigroup of bounded linear operator $(T(t))_{t\geq 0}$ defined on a Banach space $(X, ||\cdot||)$, $x_{0}\in X$, $0=t_{0}=s_{0}<t_{1}\leq s_{1} \leq t_{2} <\cdots <t_{N}\leq s_{N}\leq t_{N=1}:=a$, for all $i=1,...,N$ and $f:[0,a]\times X \rightarrow X$ is a suitable functions.

On the other hand, in 2016 Chen et al. \cite{chen}, studied the existence of mild piecewise-continuous solutions for an initial value problem of a semi-linear evolution equation with non-instantaneous impulses in the Banach space $E$, given by
\begin{equation}
\left\{ 
\begin{array}{rll}
u^{\prime }\left( t\right) +Au\left( t\right)  & = & f\left( t,u\left(
t\right) \right) ,\text{ }t\in \underset{k=0}{\overset{m}{\bigcup }}\left(
s_{k},t_{k+1}\right]  \\ 
u\left( t\right)  & = & \gamma _{k}\left( t,u\left( t\right) \right) ,\text{ 
}t\in \underset{k=1}{\overset{m}{\bigcup }}\left( t_{k},s_{k}\right]  \\ 
u\left( 0\right)  & = & u_{0}%
\end{array}%
\right. 
\end{equation}
where $A: D(A)\subset E \rightarrow E$ is a closed linear operator, $A$ is the infinitesimal generator of a strongly continuous semigroup $(T(t))_{t\geq 0}$ in $E$, $0<t_{1}<t_{2}<\cdots <t_{m}<t_{m+1}:=a$, $a>0$ is a constant, $s_{0}:=0$ and $s_{k}\in (t_{k},t_{k+1})$ for each $k=1,2,...,m$, $f:[0,a]\times E \rightarrow E$ is a given nonlinear function satisfying some assumptions, $\gamma_{k}:(t_{k},s_{k}]\times E \rightarrow E$ is non-instantaneous impulsive function for $u_{0}\in E$. 

In addition, we can also highlight the importance of the work in the context of the integrodifferential equations, which many researchers dedicate to investigate the existence and uniqueness of mild solutions of these integrodifferential equations \cite{bala,samuel,park,radha}.

With the expansion of fractional calculus \cite{RHM,KSTJ,oliveira,ZE1,ZE2}, numerous definitions of fractional derivatives were introduced and with them discussed applications and studies on the existence, uniqueness and Ulam-Hyers stability of classical and mild solutions of equations: impulsive, evolution and functional \cite{abbas2015,dabas,nonlinear,wangmitt,wang2,zhou1}. 

In 2010 Zhou and Jiao \cite{zhou} also in the context of fractional evolution equations, however neutral, presented a brief discussion on the criteria of existence and uniqueness of mild solutions. We can also highlight the important paper carried out by Shu et al. \cite{shu}, on the existence of mild solutions for a class of fractional impulsive partial differences. 
 
In 2014 Xie \cite{xie} performed out a brief study on the existence and uniqueness of mild solutions for the fractional impulse evolution integrodifferential equation with infinite decay in the Banach space, highlighting the important result on the extension of the existence theorem associated with fractional order differential equations.

In this sense, in 2015 Fu et al. \cite{fu2015}, devoted themselves to investigating the existence of mild solutions to a Cauchy problem and a fractional impulse evolution equation with non-instantaneous impulses, through the theory of semigroups and theorems of fixed points. On the other hand, results on Ulam-Hyers stabilities have been a tool of study along the slopes of several researchers, especially when it comes to mild solutions of fractional differential equations \cite{abbas2015,nonlinear,wangmitt}. For a reading of other results on existence, uniqueness and Ulam-Hyers stabilities of fractional differential equations and integrodifferential equations, we suggest \cite{chauhan,chauhan1,dabas,dabas1,wang2,zhou1}.

In this sense, it motivated us to propose a work related to the study of the existence, uniqueness and Ulam-Hyers stability of the fractional differential equation with non-instantaneous impulses given by Eq.(\ref{eq:1}), below.

Consider the fractional differential equations (FDE) with non-instantaneous impulses in a $P\beta$-normed Banach space of the form
\begin{equation}\label{eq:1}
\left\{ 
\begin{array}{cll}
^{\mathbf{H}}\mathfrak{D}_{a^{+}}^{\alpha ,\beta }u(t) & = & \mathcal{A}(t)+\mathbf{F}_{\mathfrak{T},\mathfrak{V}}(t,u(t)),\text{ }t\in \lbrack s_{i},t_{i+1}],\text{ }i\in
\lbrack 0,m] \\ 
I_{0^{+}}^{1-\gamma }u(0)+g(u) & = & u_{0} \\ 
u(t) & = & \xi (t,u(t)),t\in (t_{i},s_{i}],i\in \lbrack 1,m]%
\end{array}%
\right. 
\end{equation}
where $\mathbf{F}_{\mathfrak{T},\mathfrak{V}}(t,u(t)):=f(t,u(t),\mathfrak{T} u(t),\mathfrak{V}u(t))$, $^{\mathbf{H}}\mathfrak{D}_{a^+}^{\alpha,\beta}(\cdot)$ is the Hilfer fractional derivative of order $0<\alpha\leq{1}$ and type $0\leq{\beta}\leq{1}$, $I_{0^+}^{1-\gamma}(\cdot)$ is the Riemann-Liouville fractional integral in a Banach space $\Omega$, where $\mathcal{A}$ is the infinitesimal $\alpha$-generator of a semigroup (fractional) $\{\mathbb{P}(t)/t\geq{0}\}$ and $t_i,s_i$ are fixed numbers satisfying $0=s_0<t_1\leq{s_1}\leq{t_2}<\cdots<s_{m-1} \leq{t_m}\leq{s_m}\leq{t_{m+1}}=T$, $f:[0,T]\times \Omega\times \Omega\times \Omega\rightarrow{T}$ is continuous $g\in PC_{1-\gamma}(J,\Omega)$ and $x_i:[t_i,s_i]\times \Omega\rightarrow\Omega$ is continuous for all $i=1,2,\cdots,m$,
\begin{equation*}
\mathfrak{T}(t)=\int_{0}^{t}\mathbf{K}\left( t,s\right) u(s)ds,\,\,\mathbf{K}\in \mathbb{C}[\mathbf{D},\mathbb{R}%
^{+}],
\end{equation*}
\begin{equation*}
\mathfrak{V}(t)=\int_{0}^{t}\mathbf{H}(t,s)u(s)ds,\,\,\mathbf{H}\in \mathbb{C}[\mathbf{N},\mathbb{R}^{+}],
\end{equation*}
where $\mathbf{D}=\{(t,s)\in\mathbb{R}^{2}:0\leq{s}\leq{t}\leq{T}\}, \mathbf{N}=\{(t,s)\in\mathbb{R}^2: 0\leq{s}\leq{t}\leq{T}\}$ and $PC_{1-\gamma}(J,\Omega)$ consist of a function $u$ that are a map from $J$ into $\Omega$ such that $u(t)$ is continuous.

The main motivation to write this paper comes from the works \cite{hernandez,princ} and with the aim of proposing more general results on the existence, uniqueness and Ulam-Hyers stability in the field of fractional calculus, allowing a wider range of options and tools for future applications.

The paper is organized as follows. In Section 2 we present some spaces of functions, as well as their respective norms, as well as the definitions of $\Psi$-Riemann-Liouville fractional integral and $\Psi$-Hilfer fractional derivative. Also we present the concept of mild solution and two lemmas in the sense of integral equation and mild solution, central points of the paper. On the other hand, the definitions of stabilities: $\delta$-Ulam-Hyers, $\delta$-Ulam-Hyers generalized and $\delta$-Ulam-Hyers-Rassias are presented. To conclude the section, results on the Gronwall inequality are also presented. In section 3 we propose an extension for the Gronwall inequality in such a way make it possible to investigate $\delta$-Ulam-Hyers-Rassias stability of the fractional differential equation and in section 4 we will look at the second main purpose of this paper, the uniqueness of the $\delta$-Ulam-Hyers-Rassias stability of Eq.(\ref{eq:1}). Concluding remarks close the paper.

\section{Preliminaries} \label{sec:1}

In this section we present some spaces of functions, as well as their respective norms. In this sense, the concepts of $\Psi$-Riemann-Liouville fractional integral and $\Psi$-Hilfer fractional derivative are introduced. Since we will investigate results involving mild solution of fractional differential equations, we present the concept of mild solution and two lemmas in the sense of integral equation and mild solution. On the other hand, the definitions of stabilities: $\delta$-Ulam-Hyers, $\delta$-Ulam-Hyers generalized and $\delta$-Ulam-Hyers-Rassias are presented. To conclude the section, results on the Gronwall inequality are discussed.

Suppose $E$ is a vector over $\mathbb{R}$. A function ${\lVert \cdot \rVert}_{\delta}: E\rightarrow[0,\infty)$, $(0<\delta\leq{1})$ is called a $\delta$-norm if satisfies \cite{princ}:
\begin{enumerate}
\item ${\lVert x \rVert}_{\delta}=0$ if and only if $x=0$;
\item ${\lVert \lambda x \rVert}_{\delta}=|\lambda|^{\delta}{\lVert x \rVert}_{\delta}$, for all $\lambda\in \mathbb{K}$
and all $x\in E$;
\item ${\lVert x+y \rVert}_{\delta}\leq{\lVert x \rVert}_{\delta}+{\lVert y \rVert}_{\delta}$,
\end{enumerate}
where ${\lVert x \rVert}_{\delta}:=\displaystyle \max_{t\in{J}}|x(t)|^{\delta}$ for $x\in C(J,\Omega)$. Taking $\delta=1$, we have the space of continuous function $C(J,\Omega)$ with norm
\begin{equation*}
\left\Vert {x}\right\Vert =\max_{t\in {J}}\left\vert x(t)\right\vert .
\end{equation*}

Let $J=[0,T]$, $J'=(0,T]$ and $C(J,\Omega)$ the space of continuous function. The weighted space of functions $x$ on $J'$ is defined by \cite{ZE1}
\begin{eqnarray*}
C_{1-\gamma}(J,\Omega)=\{x\in C(J',\Omega): t^{1-\gamma}x(t)\in C(J,\Omega)\},
\end{eqnarray*}
$0\leq{\gamma}\leq{1}$, with the norm
\begin{eqnarray*}
{\lVert x \rVert}_{C_{1-\gamma};\delta}:=\sup_{t\in{J'}}|t^{1-\gamma}x(t)|^{\delta}.
\end{eqnarray*}
Obviously, the space $C_{1-\gamma}(J, \Omega)$ is a Banach space. We now introduce the definition piecewise weighted space of functions $x$ on $C_{1-\gamma}((t_k,t_{k+1}],\Omega)$ fundamental in the investigation of the main purpose of this paper, given by
\begin{equation*}
\mathcal{PC}_{1-\gamma }(J,\Omega)=\left\{ 
(t-t_{k})^{1-\gamma }x(t)\in \mathcal{C}_{1-\gamma }((t_{k},t_{k+1}],\mathbb{%
R}); \exists \lim_{t\rightarrow t_{k}}(t-t_{k})^{1-\gamma }x(t),\text{for }k=1,\ldots ,m%
\right\} 
\end{equation*}
where the norm is given by
\begin{equation*}
{\lVert x}\rVert _{\mathcal{PC}_{1-\gamma };\delta }:=\max_{k=1,\ldots
,m}\left\{ \sup_{t\in (t_{k},t_{k+1}]}\left\vert (t-t_{k})^{1-\gamma
}x(t)\right\vert ^{\delta }\right\} 
\end{equation*}
and there exists $x(t_k^-)$ and $x(t_k^+)$, $k=1,\ldots,m$ with $x(t_k^-)=x(t_k^+)$. The space $PC_{1-\gamma}(J,\Omega)$ is also a Banach space. Also, we have the following function spaces, given by \cite{ZE1}
\begin{eqnarray*}
PC_{1-\gamma}^{1}(J,\Omega)=\{x\in PC_{1-\gamma}(J,\Omega),x'\in PC_{1-\gamma}(J,\Omega)\}
\end{eqnarray*}
and
\begin{eqnarray*}
PC_{1-\gamma}^{\alpha,\beta}(J,\Omega)=\{x\in PC_{1-\gamma}(J\Omega), ^{\mathbf{H}}\mathfrak{D}_{0^+}^{\alpha,\beta}x\in PC_{1-\gamma}(J,\Omega)\}
\end{eqnarray*}
with the norm
\begin{eqnarray*}
{\lVert x \rVert}_{P{C}_{1-\gamma;\delta}^1}={\lVert x \rVert}_{PC; \delta}+{\lVert x' \rVert}_{PC_{1-\gamma};\delta}.
\end{eqnarray*}

Clearly, $PC_{1-\gamma}^{1}$ endowed with the norm ${\lVert \cdot \rVert}_{P{C}_{1-\gamma;\delta}^{1}}$ is a
$P\delta$-Banach space.

Let $J=(0,T)$ be a finite or infinite interval of the real line $\mathbb{R}_{+}$ and $\alpha>0$. Also let $\Psi(t)$ be an increasing and positive monotone function on $J^{1}=(0,T]$ having a continuous derivative $\Psi'(t)$ on $J'=(0,T)$. The left-sided fractional integral of a function $f$ with respect to a function $\Psi$ on $J=[0,T]$ is defined by \cite{ZE1}
\begin{equation*}
\mathcal{I}_{0^{+}}^{\alpha ,\Psi }u(t)=\frac{1}{\Gamma (\alpha )}%
\int_{0}^{t}\mathcal{G}\left( t,s\right) u(s)ds
\end{equation*}
where $\mathcal{G}\left( t,s\right) :=\Psi ^{\prime }\left( s\right) \left( \Psi \left(
t\right) -\Psi \left( s\right) \right) ^{\alpha -1}$.

Choosing $\Psi(t)=t$, we have the Riemann-Liouville fractional integral given by
\begin{eqnarray}
\mathcal{I}_{0^+}^{\alpha}u(t)=\frac{1}{\Gamma(\alpha)}\int_{0}^{t}(t-s)^{\alpha-1}u(s)ds.\label{RL-int}
\end{eqnarray}

We use the Eq.(\ref{RL-int}) in this paper. On the other hand, let $n-1<\alpha\leq{n}$ with $n\in\mathbb{N}$, $J=[0,T]$ the interval and $f,\Psi\in{C}^{n}([0,T],\mathbb{R})$ be two functions such that $\Psi$ is increasing and $\Psi'(t)\neq{0}$, for all $t\in{J}$. The left-sided $\Psi$-Hilfer fractional derivative $^{\mathbf{H}}\mathfrak{D}_{0^+}^{\alpha,\beta;\Psi}(\cdot)$ of a function $f$ of order $\alpha$ and type $0\leq{\beta}\leq{1}$ is defined by \cite{ZE1}
\begin{eqnarray*}
^{\mathbf{H}}\mathfrak{D}_{0^+}^{\alpha,\beta;\Psi}u(t)=\mathcal{I}_{0^+}^{\beta(n-\alpha);\Psi}\left(\frac{1}{\Psi(t)}\frac{d}{dt}\right)^{n}
{\mathcal{I}_{0^+}^{(1-\beta)(n-\alpha);\Psi}}u(t).
\end{eqnarray*}

Choosing $\Psi(t)=t$, we have the Hilfer fractional derivative, given by
\begin{eqnarray}
^{\mathbf{H}}\mathfrak{D}_{0^+}^{\alpha,\beta}u(t)=\mathcal{I}_{0^+}^{\beta(n-\alpha)}\left(\frac{d}{dt}\right)^{n}
{\mathcal{I}_{0^+}^{(1-\beta)(n-\alpha)}}u(t). \label{Hilfer-derivative}
\end{eqnarray}

In this paper, we use the so-called Hilfer fractional derivative Eq.(\ref{Hilfer-derivative}). 

In order to facilitate the development of the paper, everywhere we present the function $\psi(\cdot)$, we assume that it is a positive and increasing monotone function.

\begin{definition} \cite{lizama} Suppose that $\mathcal{A}$ is a closed, densely defined linear operator on $\Omega$. A family $\{\mathbb{\mathbb{P}}_{\alpha}(t)/t\geq{0}\}\subset B(\Omega)$ is called an $\alpha$-times resolvent family generated by $\mathcal{A}$ if the following conditions are satisfied:
\begin{enumerate}
\item $\mathbb{\mathbb{P}}_{\alpha}(\cdot)$ is strongly continuous on $\mathbb{R_{+}}$ and $\mathbb{\mathbb{P}}_{\alpha}(0)=I$;
\item $\mathbb{\mathbb{P}}_{\alpha}(t) D(\mathcal{A})\subset D(\mathcal{A})$ and $\mathcal{A}\,\mathbb{\mathbb{P}}_{\alpha}(t)x = \mathbb{\mathbb{P}}_{\alpha}\mathcal{A}x$ for all
$x\in D(\mathcal{A}), t\geq{0}$;
\item For all $x\in D(\mathcal{A})$ and $t\geq{0}$, $\mathbb{\mathbb{P}}_{\alpha}(t)x=x+\mathcal{I}_{0^+}^{\alpha}\mathbb{\mathbb{P}}_{\alpha}(t)\mathcal{A}x$.
\end{enumerate}
\end{definition}

\begin{lemma} {\rm \cite{fu2015,gu,sousa,yang}} The fractional nonlinear differential equation {\rm (FNDE)} \textnormal{Eq.(\ref{eq:1})} is equivalent to the integral equation
\begin{equation} \label{int-eq}
\left\{ 
\begin{array}{cll}
x(t) & = & \dfrac{t^{\gamma -1}}{\Gamma (\gamma )}(x_{0}-g(u))+\dfrac{1}{%
\Gamma (\alpha )}\displaystyle\int_{0}^{t}(t-s)^{\alpha
-1}(\mathcal{A}u(s)+\mathbf{F}_{\mathfrak{T},\mathfrak{V}}(s,u(s)))ds,\,\,\,\,\,\,\,\,t\in \lbrack 0,t_{1}]; \\ 
x(t) & = & \xi_{i}(t,u(t)),\,\,\,\,\,\,\,\,t\in (t_{i},s_{i}], i=1,...,m; \\ 
x(t) & = & {\xi }_{i}(t,u(t)+\dfrac{1}{\Gamma (\alpha )}%
\displaystyle\int_{0}^{t}(t-s)^{\alpha -1}(\mathcal{A}u(s)+\mathbf{F}_{\mathfrak{T},\mathfrak{V}}(s,u(s)))ds,\,\,\,\,\,t\in
(s_{i},t_{i+1}],\,\,\,i=1,\ldots ,m.%
\end{array}%
\right.  
\end{equation}

The following is the definition of the Wright function, fundamental in mild solutions of the {\rm Eq.(\ref{eq:1})}. Then, the Wright function $M_{\alpha}(Q)$ is defined by \cite{gu,yang}
\begin{eqnarray*}
M_{\alpha}(Q)=\sum_{n=1}^{\infty}\frac{(-Q)^{n-1}}{(n-1)!\Gamma(1-\alpha n)}, \quad 0<\alpha<1,\quad Q\in\mathbb{C},
\end{eqnarray*}
satisfying the equation
\begin{eqnarray*}
\int_{0}^{\infty}\theta^{\bar{\delta}}M_{\alpha}(\theta)d\theta=\frac{\Gamma(1+\bar{\delta})}{\Gamma(1+\alpha\bar{\delta})},\quad \mbox{for} \quad \theta\geq{0}.
\end{eqnarray*}
\end{lemma}

\begin{lemma} {\rm \cite{fu2015,gu,sousa,yang}} If the integral equation \textnormal{Eq.(\ref{int-eq})} holds, then we have
\begin{equation}
\left\{ 
\begin{array}{cll}
x(t) & = & \displaystyle \mathbb{\mathbb{P}}_{\alpha ,\beta }(t)(x_{0}-g(u))+\int_{0}^{t}\mathbf{K}_{\alpha }(t-s)\mathbf{F}_{\mathfrak{T},\mathfrak{V}}(s,u(s)))ds,\,\,\,\,\,\,\,\,t\in
\lbrack 0,t_{1}]; \\ 
x(t) & = & \xi_{i}(t,u(t)),\,\,\,\,\,\,\,\,t\in (t_{i},s_{i}], i=1,...,m; \\ 
x(t) & = & \displaystyle \mathbb{\mathbb{P}}_{\alpha ,\beta }(t)[\xi
_{i}(s_{i};u(s_{i})-g(u)]+\int_{s_{i}}^{t}\mathbf{K}_{\alpha
}(t-s)\mathbf{F}_{\mathfrak{T},\mathfrak{V}}(s,u(s)))ds,\,\,\,\,\,t\in (s_{i},t_{i+1}),\,\,\,i=1,\ldots ,m.%
\end{array}%
\right. 
\end{equation}
where $\mathbf{K}_{\alpha}(t)=t^{\gamma-1}\mathbb{G}_{\alpha}(t)$, $\displaystyle \mathbb{G}_{\alpha}(t)=\int_{0}^{\infty}\alpha\theta M_{\alpha}(\theta)\mathbb{\mathbb{P}}(t^{\alpha}\theta)d\theta$,
and $\mathbb{\mathbb{P}}_{\alpha,\beta}(t)=\mathcal{I}_{\theta}^{\beta(1-\alpha)}\mathbf{K}_{\alpha}(t)$.
\end{lemma}

As one of the main purpose of this paper is to investigate the stability of the $\delta$-Ulam-Hyers-Rassias type of Eq. (\ref{eq:1}), we need some definitions in this context. For the following definitions that will be presented next, we use the paper to adapt in the fractional sense \cite{princ}.

Let $0<\delta<1$, $\epsilon>0$, $\phi\geq{0}$ and $PC_{1-\gamma}(J,\mathbb{R}_{+})$. We consider the following inequalities:
\begin{equation}\label{ineq}
\left. 
\begin{array}{c}
\left\vert {^{\mathbf{H}}\mathfrak{D}_{0^{+}}^{\alpha ,\beta }}%
v(t)-\mathcal{A}v(t)-\mathbf{F}_{\mathfrak{T},\mathfrak{V}}(s,v(s))\right\vert \leq {\epsilon },\,\,t\in
(s_{i},t_{i+1}],\,\,i\in \lbrack 0,m] \\ 
\left\vert v(t)-\xi _{i}(t,v(t))+g(t)\right\vert \leq {\epsilon },\quad t\in
(t_{i},s_{i}],\quad i\in \lbrack 1,m]%
\end{array}%
\right.   
\end{equation}%
and 
\begin{equation}\label{ineq-2}
\left. 
\begin{array}{c}
\left\vert {^{\mathbf{H}}\mathfrak{D}_{0^{+}}^{\alpha ,\beta }}%
v(t)-\mathcal{A}v(t)-\mathbf{F}_{\mathfrak{T},\mathfrak{V}}(s,v(s))\right\vert \leq {\varphi (t)},\,\,t\in
(s_{i},t_{i+1}],\,\,i\in \lbrack 0,m] \\ 
\left\vert v(t)-\xi _{i}(t,v(t))+g(t)\right\vert \leq {\phi },\quad t\in
(t_{i},s_{i}],\quad i\in \lbrack 1,m]%
\end{array}%
\right.   
\end{equation}%
and
\begin{equation}\label{ineq-3}
\left. 
\begin{array}{c}
\left\vert {^{\mathbf{H}}\mathfrak{D}_{0^{+}}^{\alpha ,\beta }}%
v(t)-\mathcal{A}v(t)-\mathbf{F}_{\mathfrak{T},\mathfrak{V}}(s,v(s))\right\vert \leq {\epsilon \varphi (t)}%
,\,\,t\in (s_{i},t_{i+1}],\,\,i\in \lbrack 0,m] \\ 
\left\vert v(t)-\xi _{i}(t,v(t))+g(t)\right\vert \leq {\epsilon \phi },\quad
t\in (t_{i},s_{i}],\quad i\in \lbrack 1,m].%
\end{array}%
\right.   
\end{equation}

\begin{definition} \cite{princ} The \textnormal{Eq.(\ref{eq:1})} is $\delta$-Ulam-Hyers stable if there exists a real number $c_{f_{j},\bar{M},\bar{G},\delta,\xi_{i},\varphi}>0$ such that for each $\epsilon>0$ and for each solution $y\in{P}C_{1-\gamma}^{1}(J,\Omega)$ of \textnormal{Eq.(\ref{ineq})} there exists a solution $x\in{P}C_{1-\gamma}^{1}(J,\Omega)$ of \textnormal{Eq.(\ref{eq:1})} with
\begin{equation*}
\left\vert v(t)-u(t)\right\vert ^{\delta }\leq {c_{f_{j},\bar{M},\bar{G}%
,\delta ,\xi _{i},\varphi }}(\epsilon ^{\delta }),\quad t\in {J}.
\end{equation*}
\end{definition}

\begin{definition} \cite{princ} The \textnormal{Eq.(\ref{eq:1})} is generalized $\delta$-Ulam-Hyers stable if there exists $\theta_{f_{j},\bar{M},\bar{G},\delta,\xi_{i},\varphi}\in C_{1-\gamma}(\mathbb{R}_{+},\mathbb{R}_{+})$, $\theta_{f_{j},\bar{M},\bar{G},\delta,\xi_{i},\varphi}(0)=0$ such that for each solution $y\in PC_{1-\gamma}^{1}(J,\Omega)$ of \textnormal{Eq.(\ref{ineq})} there exists a solution $x\in{PC}_{1-\gamma}^{1}(J,\Omega)$
of \textnormal{Eq.(\ref{eq:1})} with 
\begin{eqnarray*}
|v(t)-u(t)|^{\delta}\leq{\theta_{f_{j},\bar{M},\bar{G},\delta,\xi_{i},\varphi}}(\epsilon^{\delta}), \quad t\in{J}.
\end{eqnarray*}
\end{definition}

\begin{definition} \cite{princ} The \textnormal{Eq.(\ref{eq:1})} is $\delta$-Ulam-Hyers-Rassias stable with respect to $(\varphi,\phi)$ if there exists $c_{f_{j},\bar{M},\bar{G},\delta,\xi_{i},\varphi}>0$ such that for each $\epsilon>0$ and for each solution $y\in PC_{1-\gamma}^{1}(J,\Omega)$ of \textnormal{Eq.(\ref{ineq-3})} there exists a solution $x\in PC_{1-\gamma}^{1}(J,\Omega)$ of \textnormal{Eq.(\ref{eq:1})} with 
\begin{eqnarray*}
|v(t)-u(t)|^{\delta}\leq{c_{f_{j},\bar{M},\bar{G},\delta,\xi_{i},\varphi}}\epsilon^{\delta}(\phi^{\delta}-\varphi^{\delta}(t)), \quad t\in{J}.
\end{eqnarray*}
\end{definition}

\begin{definition} \cite{princ} The \textnormal{Eq.(\ref{eq:1})} is generalized $\delta$-Ulam-Hyers-Rassias stable with respect to $(\varphi,\phi)$ if there exists $c_{f_{j},\bar{M},\bar{G},\delta,\xi_{i},\varphi}>0$ such that for each solution $y\in PC_{1-\gamma}^{1}(J,\Omega)$ of \textnormal{Eq.(\ref{ineq-2})} there exists a solution $x\in PC_{1-\gamma}^{1}(J,\Omega)$ of \textnormal{Eq.(\ref{eq:1})} with 
\begin{eqnarray*}
|v(t)-u(t)|^{\delta}\leq{c_{f_{j},\bar{M},\bar{G},\delta,\xi_{i},\varphi}}(\phi^{\delta}+\varphi^{\delta}(t)), \quad t\in{J}.
\end{eqnarray*}
\end{definition}

\textbf{Remark 1.} A function $v\in PC_{1-\gamma}^{1}(J,\mathbb{R})$ is a solution of \textnormal{Eq.(\ref{ineq-2})} if, and only if, there is $H\in PC_{1-\gamma}(J,\mathbb{R})$ and a sequence $H_i$, $i=1,2,\ldots,m$ (which depend on $y$) such that
\begin{enumerate}
\item $|H(t)|\leq{\varphi(t)}$, $t\in{J}$ and $|H_i|\leq\phi$, $i=1,2,\ldots,m$;

\item ${^{\mathbf{H}}\mathfrak{D}_{0^+}^{\alpha,\beta}}v(t)=\mathcal{A}(t)+f(t,v(t))+G(t)$, $t\in(s_i,t_{i+1}]$, $i=1,2,\ldots,m$;

\item $v(t)=g_{i}(t,v(t))+H_i$, $t\in(t_i,s_i]$, $i=1,2,\ldots,m$.
\end{enumerate}

\textbf{Remark 2.} If $v\in PC_{1-\gamma}^{1}(J,\Omega)$ is a solution of \textnormal{Eq.(\ref{ineq-2})} then $v$ satisfies the inequality
\begin{eqnarray}\label{ineq-rk}
&&\left\vert v(t)-\mathbb{\mathbb{P}}_{\alpha ,\beta }(t-s)(\xi _{i}(t,v(t))-g(t))\right\vert
\leq {M_{1}}\phi ,\quad t\in (t_{i},s_{i}],\quad i=[1,m],  \notag \\
&&\left\vert v(t)-\mathbb{\mathbb{P}}_{\alpha ,\beta }(t-s)(v(0)-y(t))-\int_{0}^{t}\mathbf{K}_{\alpha
}(t-s)\mathbf{F}_{\mathfrak{T},\mathfrak{V}}(s,v(s))ds\right\vert    \\
&\leq &\int_{0}^{t}{\lVert \mathbf{K}_{\alpha }(t-s)\rVert }\varphi (s)ds,\quad t\in
\lbrack 0,t_{1}],  \notag \\
&&\left\vert v(t)-\mathbb{\mathbb{P}}_{\alpha ,\beta }(t-s)(\xi
_{i}(t,v(t))-g(t))-\int_{s_{i}}^{t}\mathbf{K}_{\alpha
}(t-s)\mathbf{F}_{\mathfrak{T},\mathfrak{V}}(s,v(s))ds\right\vert   \notag \\
&\leq &{\lVert \mathbb{\mathbb{P}}_{\alpha ,\beta }(t-s)\rVert }_{PC_{1-\gamma }}\phi
+\int_{s_{i}}^{t}{\lVert \mathbf{K}_{\alpha }(t-s)\rVert }_{PC_{1-\gamma }}\varphi
(s)ds,\quad t\in \lbrack s_{i},t_{i+1}],\quad i\in \lbrack 1,m].  \notag
\end{eqnarray}

In order to investigate the $\delta$-Ulam-Hyers-Rassias stability, the Gronwall inequality, as presented in Theorem \ref{theo-1} below is of paramount importance. Lemma \ref{lemma3} and Lemma \ref{lemma4} are also of great importance, since they are direct consequences of Theorem \ref{theo-1}.

\begin{teorema}\label{theo-1} {\rm\cite{gron}} Let $u,v$ be two integrable functions and $g$ continuous with domain $I:=[a,T]$. Let $\Psi\in{C}^{1}(I)$ an increasing function such that $\Psi'(t)\neq{0}$, $t\in{I}$. Assume that
\begin{enumerate}
\item $u$ and $v$ are nonnegative;
\item $g$ is nonnegative and nondecreasing.
\end{enumerate}

If
\begin{eqnarray*}
u(t)\leq{v(t)}+g(t)\int_{a}^{t}\mathcal{G}\left( t,s\right)u(s)ds,
\end{eqnarray*}
then
\begin{eqnarray*}
u(t)\leq{v(t)}+\int_{a}^{t}\sum_{k=1}^{\infty}\frac{[g(t)\Gamma(\alpha)]^k}
{\Gamma(\alpha k)}\Psi'(s)(\Psi(t)-\Psi(s))^{\alpha k -1}ds.
\end{eqnarray*}
\end{teorema}

\begin{lemma} \label{lemma3} {\rm\cite{gron}} Under the hypothesis of Theorem \ref{theo-1}, let $v$ be a nondecreasing function on $I$. Then, we have
\begin{eqnarray*}
u(t)\leq{v(t)}\,\mathbb{E}_{\alpha}[\widetilde{\Psi}^{\alpha}_{g}(t,a)], \quad \forall t\in{I},
\end{eqnarray*}
where $\widetilde{\Psi}^{\alpha}_{g}(t,a):=g(t)\Gamma(\alpha)(\psi(t)-\psi(a))^{\alpha}$, and $\mathbb{E}_{\alpha}(\cdot)$ is the classical Mittag-Leffler function.
\end{lemma}
\newpage

\begin{lemma} \label{lemma4} {\rm\cite{nonlinear}}
\begin{enumerate}
\item The following inequality holds
\begin{eqnarray*}
u(t)\leq a(t)+g(t)\int_{0}^{t}(t-s)^{\alpha-1}u(s)ds, \quad t\geq{0}
\end{eqnarray*}
where $u,a\in\mathcal{PC}(\mathbb{R}_{+},\mathbb{R}_{+})$, $a$ is nondecreasing, $u$ and $a$ nonnegative. Then, for $t\in\mathbb{R}_{+}$; we have
\begin{eqnarray*}
u(t)\leq a(t)\,\mathbb{E}_{\alpha}[g(t)\Gamma(\alpha)(t)^{\alpha}].
\end{eqnarray*}
\item Assume
\begin{eqnarray*}
u(t)\leq a(t)+\delta\,u(t)+g(t)\int_{0}^{t}(t-s)^{\alpha-1}u(s)ds+
\sum_{0<t_k<t}\beta_{k}\,u({t_k}^{-}), \quad t\geq{0}
\end{eqnarray*}
where $u,\delta,a,b\in PC(\mathbb{R}_{+},\mathbb{R}_{+})$, $a,\delta$ is
nondecreasing, $\beta_k>0$, $k=1,\ldots,m$ and $u$ and $a$ nonnegative. Then, for $t\in\mathbb{R}_{+}$, we have
\begin{eqnarray*}
u(t)\leq{a(t)}[\delta\,\mathbb{E}_{\alpha}[g(t)\Gamma(\alpha)t^{\alpha}]
(1+\beta\,\mathbb{E}_{\alpha}[g(t)\Gamma(\alpha)t^{\alpha}])^{k-1}+
(1+\beta\,\mathbb{E}_{\alpha}[g(t)\Gamma(\alpha)t^{\alpha}])^{k}]
\mathbb{E}_{\alpha}[g(t)\Gamma(\alpha)t^{\alpha}].
\end{eqnarray*}
\end{enumerate}
\end{lemma}

\section{Extension of the Gronwall inequality}

The Gronwall inequality is an important tool in the research of data dependence and Ulam-Hyers stabilities. When proposing more general fractional derivatives, one of the next objectives is to propose an extension for the Gronwall inequality in such a way to make it possible to investigate Ulam-Hyers stabilities of certain fractional differential equations. In this sense, this section is intended to propose an extension for the Gronwall inequality and through this, to present some particular cases in the form of lemmas, that will be used in section 4.

\begin{teorema}
Let $u\in PC_{1-\gamma,{\Psi}}(I,\mathbb{R}_{+})$ satisfying the following inequality
\begin{equation*}
u(t)\leq v(t)+\delta u(t)+g(t)\int_{a}^{t}\mathcal{G}\left( t,s\right)u(s)ds+\sum_{a<t_{k}<t}\beta _{k}\,u(t_{k}^{-}),\,\,t\geq {a},
\end{equation*}
where $\mathcal{G}\left( t,s\right) :=\Psi ^{\prime }\left( s\right) \left( \Psi \left(
t\right) -\Psi \left( s\right) \right) ^{\alpha -1}$ and $v(t), u(t), \delta$ are nonnegative and $\beta_k>0$, $k=1,\dots,m$ nonnegative constants too. Then
\begin{eqnarray}
u(t) &\leq &{v(t)}\left[ \delta \,\mathbb{E}_{\alpha }[\widetilde{\Psi}^{\alpha}_{g}(t,a)]\prod_{i=1}^{k-1}\left( 1+\beta _{i}\,\mathbb{E}_{\alpha }[\widetilde{\Psi}^{\alpha}_{g}(t_{i},a)]\right)
\right. \left. +\prod_{i=1}^{k}\left( 1+\beta _{i}\,\mathbb{E}_{\alpha
}[\widetilde{\Psi}^{\alpha}_{g}(t_{i},a)]\right) \right]
\mathbb{E}_{\alpha }[\widetilde{\Psi}^{\alpha}_{g}(t,a)], \notag \\  \label{eq0-th2}
\end{eqnarray}
with $t\in(t_k,t_{k+1}]$ and $\widetilde{\Psi}^{\alpha}_{g}(t_{i},a):=g(t)\Gamma(\alpha)(\psi(t_{i})-\psi(a))^{\alpha}$ with $i=1,2$.

\begin{proof} Consider $t\in I=[a,T]$ and by means of the {\rm Lemma \ref{lemma4}}, we get
\begin{eqnarray}
u(t)\leq{[v(t)+\delta u(t)]}\,\mathbb{E}_{\alpha}[\widetilde{\Psi}^{\alpha}_{g}(t,a)], \quad t\in[a,t_1]
\label{eq1-th2}
\end{eqnarray}
and
\begin{eqnarray}
u(t)\leq\left(v(t)+\delta u(t)+\sum_{i=1}^{k}\beta_{i}u(t_{i}^{+}\right)\,
\mathbb{E}_{\alpha}[\widetilde{\Psi}^{\alpha}_{g}(t,a)], \quad t\in(t_{k},t_{k+1}],
\nonumber\\
\label{eq2-th2}
\end{eqnarray}
where $k=1,2,\ldots,m$. On the other hand, let $t\in(t_1,t_2)$ and using the \textnormal{Eq.(\ref{eq1-th2})}, we have
\begin{eqnarray}
u(t)&\leq &[v(t)+\delta u(t)+\beta_{1}u(t_{1}^{-})]\,\mathbb{E}_{\alpha}[\widetilde{\Psi}^{\alpha}_{g}(t,a)]
\nonumber\\
&\leq &\biggl( v(t)+\delta v(t)\,\mathbb{E}_{\alpha}[\widetilde{\Psi}^{\alpha}_{g}(t,a)]
+\beta_{1}v(t)\,\mathbb{E}_{\alpha}[g(t)\Gamma(\alpha)(\Psi(t_1)-\Psi(a))^{\alpha}]\biggr) \nonumber \,\mathbb{E}_{\alpha}[\widetilde{\Psi}^{\alpha}_{g}(t,a)]\nonumber\\
&=&v(t)\biggl\{1+\delta\,\mathbb{E}_{\alpha}[\widetilde{\Psi}^{\alpha}_{g}(t,a)]+
\beta_{1}\,\mathbb{E}_{\alpha}[\widetilde{\Psi}^{\alpha}_{g}(t_{1},a)]\biggr\} \,\mathbb{E}_{\alpha}[\widetilde{\Psi}^{\alpha}_{g}(t,a)]. \label{eq3-th2}
\end{eqnarray}

Now, for $t\in(t_2,t_3]$ and using the \textnormal{Eq.(\ref{eq1-th2})} and \textnormal{Eq.(\ref{eq3-th2})}, we get
\begin{eqnarray}
u(t)&\leq &\biggl(v(t)+\delta u(t)+\beta_{1}u(t_{1}^{-})+\beta_{2}u(t_{2}^{-})\biggr)
\,\mathbb{E}_{\alpha}[\widetilde{\Psi}^{\alpha}_{g}(t,a)]\nonumber\\
&\leq& \biggl(v(t)+\delta v(t)\,\mathbb{E}_{\alpha}[\widetilde{\Psi}^{\alpha}_{g}(t,a)]
+\beta_{1}v(t)\,\mathbb{E}_{\alpha}[\widetilde{\Psi}^{\alpha}_{g}(t_{1},a)]\nonumber\\
&+&\beta_{2}v(t)\biggl\{1+\delta\,\mathbb{E}_{\alpha}[\widetilde{\Psi}^{\alpha}_{g}(t,a)]+
\beta_{1}\,\mathbb{E}_{\alpha}[\widetilde{\Psi}^{\alpha}_{g}(t_{1},a)]\biggr\}\,\mathbb{E}_{\alpha}[\widetilde{\Psi}^{\alpha}_{g}(t_{2},a)]\biggr)
\mathbb{E}_{\alpha}[\widetilde{\Psi}^{\alpha}_{g}(t,a)]\nonumber\\
&=&v(t)\biggr\{\biggl(1+\delta\,\mathbb{E}_{\alpha}[\widetilde{\Psi}^{\alpha}_{g}(t,a)]
+\beta_{1}\mathbb{E}_{\alpha}[\widetilde{\Psi}^{\alpha}_{g}(t_{1},a)]\biggr)\biggl(1+\beta_{2}\mathbb{E}_{\alpha}[\widetilde{\Psi}^{\alpha}_{g}(t_{2},a)]\biggr)\biggr\}
\mathbb{E}_{\alpha}[\widetilde{\Psi}^{\alpha}_{g}(t,a)]. \notag \\ \label{eq4-th2}
\end{eqnarray}

Now, for $t\in(t_3,t_4]$ and using \textnormal{Eq.(\ref{eq1-th2})} and \textnormal{Eq.(\ref{eq4-th2})}, we get
\begin{eqnarray*}
u(t) &\leq &\biggl(v(t)+\delta u(t)+\beta _{1}u(t_{1}^{-})+\beta _{2}u(t_{2}^{-})+\beta _{3}u(t_{3}^{-})\biggr)\mathbb{E}_{\alpha }[\widetilde{\Psi}^{\alpha}_{g}(t,a)] \\
&\leq &\biggl(v(t)+\delta v(t)\mathbb{E}_{\alpha }[\widetilde{\Psi}^{\alpha}_{g}(t,a)]+\beta _{1}v(t)\mathbb{E}_{\alpha }[\widetilde{\Psi}^{\alpha}_{g}(t_{1},a)] + \beta _{2}v(t)\biggl\{1+\delta \,\mathbb{E}_{\alpha }[\widetilde{\Psi}^{\alpha}_{g}(t,a)]+\beta _{1}\mathbb{E}_{\alpha }[\widetilde{\Psi}^{\alpha}_{g}(t_{1},a)]\biggr\} \\
&\times &\mathbb{E}_{\alpha }[\widetilde{\Psi}^{\alpha}_{g}(t_{2},a)]+\beta _{3}v(t)\biggl\{1+\delta \,\mathbb{E}_{\alpha }[\widetilde{\Psi}^{\alpha}_{g}(t,a)]+\beta _{1}\mathbb{E}_{\alpha }[\widetilde{\Psi}^{\alpha}_{g}(t_{1},a)]+\beta _{2}\biggl(1+\delta \,\mathbb{E}_{\alpha }[\widetilde{\Psi}^{\alpha}_{g}(t,a)] \\
&+&\beta _{1}\mathbb{E}_{\alpha }[\widetilde{\Psi}^{\alpha}_{g}(t_{1},a)]\biggr)\mathbb{E}_{\alpha }[\widetilde{\Psi}^{\alpha}_{g}(t_{2},a)]\biggr\}\mathbb{E}_{\alpha }[\widetilde{\Psi}^{\alpha}_{g}(t_{3},a)]\biggr)\mathbb{E}_{\alpha }[\widetilde{\Psi}^{\alpha}_{g}(t,a)] \\
\end{eqnarray*}

\begin{eqnarray*}
&=&v(t)\biggl\{\biggl(1+\delta \,\mathbb{E}_{\alpha }[\widetilde{\Psi}^{\alpha}_{g}(t,a)]+\beta _{1}\mathbb{E}_{\alpha }[\widetilde{\Psi}^{\alpha}_{g}(t_{1},a)] +\beta _{2}\biggl(1+\delta \,\mathbb{E}_{\alpha }[\widetilde{\Psi}^{\alpha}_{g}(t,a)]+\beta _{1}\mathbb{E}_{\alpha }[\widetilde{\Psi}^{\alpha}_{g}(t_{1},a)]\biggr) \\
&\times &\mathbb{E}_{\alpha }[\widetilde{\Psi}^{\alpha}_{g}(t_{2},a)]+\beta _{3}\biggl[1+\delta \,\mathbb{E}_{\alpha }[\widetilde{\Psi}^{\alpha}_{g}(t,a)] +\beta _{1}\mathbb{E}_{\alpha }[\widetilde{\Psi}^{\alpha}_{g}(t_{1},a)]+\beta _{2}\biggl(1+\delta \,\mathbb{E}_{\alpha }[\widetilde{\Psi}^{\alpha}_{g}(t,a)] \\
&+&\beta _{1}\mathbb{E}_{\alpha }[g(t)\Gamma (\alpha )(\Psi (t_{1})-\Psi
(a))^{\alpha }]\biggr)\mathbb{E}_{\alpha }[\widetilde{\Psi}^{\alpha}_{g}(t_{2},a)]\biggr]\mathbb{E}_{\alpha }[\widetilde{\Psi}^{\alpha}_{g}(t_{3},a)]\biggr\} \mathbb{E}_{\alpha }[\widetilde{\Psi}^{\alpha}_{g}(t,a)]\\
&=&v(t)\biggl[\delta \,\mathbb{E}_{\alpha }[\widetilde{\Psi}^{\alpha}_{g}(t,a)]\biggl(1+\beta _{2}\mathbb{E}_{\alpha }[\widetilde{\Psi}^{\alpha}_{g}(t_{2},a)]\biggr) \biggl(1+\beta _{3}\mathbb{E}_{\alpha }[\widetilde{\Psi}^{\alpha}_{g}(t_{3},a)]\biggr)+\biggl(1+\beta _{1}\mathbb{E}_{\alpha }[\widetilde{\Psi}^{\alpha}_{g}(t,a)]\biggr) \\
&\times &\biggl(1+\beta _{2}\mathbb{E}_{\alpha }[\widetilde{\Psi}^{\alpha}_{g}(t_{2},a)]\biggr)\biggl(1+\beta _{3}\mathbb{E}_{\alpha }[\widetilde{\Psi}^{\alpha}_{g}(t_{3},a)]\biggr)\biggr] \mathbb{E}_{\alpha }[\widetilde{\Psi}^{\alpha}_{g}(t,a)] \\
&=&v(t)\biggl[\delta \mathbb{E}_{\alpha }[\widetilde{\Psi}^{\alpha}_{g}(t,a)]\prod_{i=1}^{2}\biggl(1+\beta _{i}\mathbb{E}_{\alpha }[\widetilde{\Psi}^{\alpha}_{g}(t_{i},a)] +\prod_{i=1}^{3}\biggl(1+\beta _{i}\mathbb{E}_{\alpha }[\widetilde{\Psi}^{\alpha}_{g}(t_{i},a)]\biggr)\biggr]\mathbb{E}_{\alpha }[\widetilde{\Psi}^{\alpha}_{g}(t,a)],\quad t\in (t_{3},t_{4}].
\end{eqnarray*}

Following the same procedure, we have
\begin{eqnarray*}
u(t)&\leq& v(t)\biggl[\delta\,\mathbb{E}_{\alpha}[\widetilde{\Psi}^{\alpha}_{g}(t,a)]\prod_{i=1}^{n-1}\biggl(1+\beta_{i}
\mathbb{E}_{\alpha}[\widetilde{\Psi}^{\alpha}_{g}(t_{i},a)]\biggr)+\prod_{i=1}^{n}\biggl(1+\beta_{i}\mathbb{E}_{\alpha}[g\widetilde{\Psi}^{\alpha}_{g}(t_{i},a)]\biggr)\bigg]
\mathbb{E}_{\alpha}[\widetilde{\Psi}^{\alpha}_{g}(t,a)].
\end{eqnarray*}
\end{proof}
\end{teorema}

\begin{corolario}\label{Cor-A} Let $u\in PC_{1-\gamma,{\Psi}}(I,\mathbb{R}_{+})$ satisfying the following inequality
\begin{eqnarray*}
u(t)\leq c+\delta\,u(t)+g(t)\int_{a}^{t}\mathcal{G}\left( t,s\right)u(s)ds+\sum_{a<t_k<t}\beta_{k}u(t_{k}^{-}), \quad t\geq{a},
\end{eqnarray*}
where $\mathcal{G}\left( t,s\right) :=\Psi ^{\prime }\left( s\right) \left( \Psi \left(
t\right) -\Psi \left( s\right) \right) ^{\alpha -1}$ and $\delta,u(t)\in PC_{1-\gamma,\Psi}(I,\mathbb{R}_{+})$ nonnegative and $c,\beta_{k}>0$, $k=1,2,\ldots,m$ are nonnegative constants. Then
\begin{eqnarray*}
u(t)&\leq&{c}\biggl[\delta\,\mathbb{E}_{\alpha}[\widetilde{\Psi}^{\alpha}_{g}(t,a)]\prod_{i=1}^{k-1}\biggr(1+\beta_{i}
\mathbb{E}_{\alpha}[\widetilde{\Psi}^{\alpha}_{g}(t_{i},a)]\biggr)=\prod_{i=1}^{k}\biggl(1+\beta_{i}\mathbb{E}_{\alpha}[\widetilde{\Psi}^{\alpha}_{g}(t_{i},a)]\biggr)\biggr]
\mathbb{E}_{\alpha}[\widetilde{\Psi}^{\alpha}_{g}(t,a)], \quad (t_k,t_{k+1}].
\end{eqnarray*}
\end{corolario}
\begin{proof}
The proof follows the same steps of Theorem \ref{theo-1}.
\end{proof}

\begin{corolario}Let $u\in PC_{1-\gamma}(I,\mathbb{R}_{+})$ satisfying the following inequality
\begin{eqnarray}
u(t)\leq{v(t)}+g(t)\int_{a}^{t}\mathcal{G}\left( t,s\right)u(s)ds+\sum_{a<t_k<t}\beta_{k}u(t_k), \quad t\geq{a},
\end{eqnarray}
where $\mathcal{G}\left( t,s\right) :=\Psi ^{\prime }\left( s\right) \left( \Psi \left(
t\right) -\Psi \left( s\right) \right) ^{\alpha -1}$ and $v(t)\in PC_{1-\gamma,\Psi}(I,\mathbb{R}_{+})$ are nonnegative and $\beta_k>0$, $k=1,\ldots,m$ nonnegative constant too. Then
\begin{eqnarray*}
u(t)\leq{v(t)}\biggl[\prod_{i=1}^{k}\biggr(1+\beta_{i}\mathbb{E}_{\alpha}[\widetilde{\Psi}^{\alpha}_{g}(t_{i},a)]\biggr)\biggr]
\mathbb{E}_{\alpha}[\widetilde{\Psi}^{\alpha}_{g}(t,a)], \quad t\in(t_k,t_{k+1}].
\end{eqnarray*}
\end{corolario}
\begin{proof}
The proof follows the same steps of Theorem \ref{theo-1}.
\end{proof}

\begin{corolario}Let $u\in PC_{1-\gamma}(I,\mathbb{R}_{+})$ satisfying the following inequality
\begin{eqnarray*}
u(t)\leq{v(t)}+\delta u(t)+g(t)\int_{a}^{t}(t-s)^{\alpha-1}u(s)ds+\sum_{a<t_k<t}\beta_{k}u(t_{k}^{-}), \quad t\geq{a},
\end{eqnarray*}
where $v(t),u(t),\delta\in PC_{1-\gamma}(I,\mathbb{R}_{+})$ are nonnegative and $\beta_{k}>0$, $k=1,\ldots,m$, too nonnegative constants. Then
\begin{eqnarray*}
u(t)&\leq&{v(t)}\bigg[\delta\,\mathbb{E}_{\alpha}[g(t)\Gamma(\alpha)(t-a)^{\alpha}]\prod_{i=1}^{k-1}\bigg(1+\beta_{i}
\mathbb{E}_{\alpha}[g(t)\Gamma(\alpha)(t_i-a)^{\alpha}]\bigg)\\
&+&\prod_{i=1}^{k}\bigg(1+\beta_{i}\mathbb{E}_{\alpha}[g(t)\Gamma(\alpha)(t_i-a)^{\alpha}]\bigg)\bigg]
\mathbb{E}_{\alpha}[g(t)\Gamma(\alpha)(t-a)^{\alpha}], \quad t\in(t_k,t_{k+1}].
\end{eqnarray*}
\end{corolario}
\begin{proof}
The proof follows the same steps of Theorem \ref{theo-1}.
\end{proof}

Here we present some results that are particular cases of Theorem \ref{theo-1}.

\noindent \textbf{Remark 1.} \begin{enumerate} \item Taking $\delta=0$, $g(t)=\lambda$ and $\Psi(t)=t$ in the 
\textnormal{Eq.(\ref{eq0-th2})}, we have the Corollary 2.7 \cite{choi2015};

\item Taking $\delta=0$, $v(t)=c$, $g(t)=\lambda$ and $\Psi(t)=t$ in the \textnormal{Eq.(\ref{eq0-th2})}, we have the Theorem 2.6 \cite{choi2015};

\item Taking $\Psi(t)=\ln t$ in the \textnormal{Eq.(\ref{eq0-th2})}, 
we have the following

\begin{eqnarray*}
u(t)&\leq&{v(t)}\biggl[\delta\,\mathbb{E}_{\alpha}[g(t)\Gamma(\alpha)(\ln t-\ln a)^{\alpha}]\prod_{i=1}^{k-1}
\biggl(1+\beta_{i}\mathbb{E}_{\alpha}[g(t)\Gamma(\alpha)(\ln t_{i}-\ln a)^{\alpha}]\biggr)\\
&+&\prod_{i=1}^{k}\biggl(1+\beta_{i}\mathbb{E}_{\alpha}[g(t)\Gamma(\alpha)(\ln t-\ln a)^{\alpha}]\biggr)\biggr]
\mathbb{E}_{\alpha}[g(t)\Gamma(\alpha)(\ln t-\ln a)^{\alpha}], \quad t\in(t_k,t_{k+1}].
\end{eqnarray*}

\item In the recent paper Sousa and Oliveira \cite{ZE1}, introducing a new fractional derivative and presenting a part of the vast class of the fractional derivatives derived from it, also presented some particular cases of the $\Psi$-Riemann-Liouville fractional integral. In this sense, from the choice of the $\Psi(t), v(t)$ and $\delta$, it is possible to obtain other types of inequalities from Theorem \ref{theo-1} and which can be useful for another open questions.
\end{enumerate}
\section{Uniqueness and $\delta$-Ulam-Hyers-Rassias stability}

In this section we will investigate the second objective of this paper, the uniqueness of the $\delta$-Ulam-Hyers-Rassias stability of Eq.(\ref{eq:1}).

To get our key results, we need the following conditions:

\noindent (H1): $\mathcal{A}$ is the infinitesimal generator of a strongly continuous semigroup fractional $\mathbb{\mathbb{P}}_{\alpha,\beta}(t)$, whose domain $D(\mathcal{A})$ is dense in $H$ such that ${\lVert \mathbb{\mathbb{P}}_{\alpha,\beta}(t) \rVert}_{PC_{1-\gamma}}\leq{\mathfrak{M}}$, for all $t\in{J}$.\\

\noindent (H2): $f\in C_{1-\gamma}(J\times \Omega \times \Omega\times \Omega\rightarrow{\Omega})$,
$g:\Omega\rightarrow\Omega$ and there exists constants ${L}_{f_1}, {L}_{f_2}, {L}_{f_3}\geq{0}$,
$\tilde{{L}}\geq{0}$ such that
\begin{eqnarray*}
&|f(t,x_1,x_2,x_3)-f(t,y_1,y_2,y_3)|\leq{{L}_{f_1}}|x_1-y_1|+{L}_{f_2}|x_2-y_2|+{L}_{f_3}|x_3-y_3|,\\
&t\in{J}, \quad x_j,y_j\in\Omega, \,\,j=1,2,3.\\
&| g(\eta_1)-g(\eta_2)|\leq{\tilde{L}}|\eta_1-\eta_2|,\quad \eta_1,\eta_2\in PC_{1-\gamma}(J,\Omega)
\end{eqnarray*}

\noindent (H3): Denote $M_f=\max\{M_{f_1},M_{f_2},M_{f_3}\}$,
\begin{eqnarray*}
\mathfrak{F}_{1}^{*}&=&\sup_{t\in{J}}\int_{0}^{t}(t-s)^{1-\alpha}|\mathfrak{F}_{1}(t,s)|ds\leq{\infty},\\
\mathfrak{F}_{2}^{*}&=&\sup_{t\in{J}}\int_{0}^{t}(t-s)^{1-\alpha}|\mathfrak{F}_{2}(t,s)|ds\leq{\infty},\\
\mathfrak{F}_{3}^{*}&=&\sup_{t\in{J}}\int_{0}^{t}(t-s)^{1-\alpha}ds.
\end{eqnarray*}

\noindent (H4): $\xi_{i}\in C_{1-\gamma}([t_i,s_i]\times\Omega,\Omega)$ and there are positive $L_{\xi_i}>$0, $i=1,2,\ldots,m$ such that
\begin{eqnarray*}
|\xi_{i}(t,u_1)-\xi_{i}(t,u_2)|\leq{L_{\xi_i}}|u_1-u_2|
\end{eqnarray*}

\noindent for each $t\in[t_i,s_i]$ and all $u_1,u_2\in\mathbb{R}$.

\noindent (H5): Let $\varphi\in C_{1-\gamma}(J,\mathbb{R}_{+})$ be a nondecreasing function. There exits
$c_{\varphi}>0$ such that 
\begin{eqnarray*}
\int_{0}^{t}\varphi(s)ds\leq c_{\varphi}\,\varphi(t),
\end{eqnarray*}
for each $t\in{J}$.

\begin{teorema} \label{th3} Assume that \textnormal{(H1)-(H2)} are satisfied. Then \textnormal{Eq.(\ref{int-eq})} has a unique solution provided that
\begin{eqnarray}
\Lambda&:=&max\bigg\{\mathfrak{M}({L}_{\xi_i}^{\delta}+\tilde{L}^{\delta})+\mathfrak{M}
\bigg(L_{f_1}^{\delta}(t_{i+1}-s_i)^{\delta}+L_{f_2}^{\delta}\frac{(t_{i+1}-s_i)^{\alpha\delta}}{\alpha}\,\mathfrak{F}_{1}^{*}
+L_{f_3}^{\delta}\frac{(t_{i+1}-s_i)^{\alpha\delta}}{\alpha}\,\mathfrak{F}_{2}^{*}\bigg);\nonumber\\
&&\mathfrak{M}\bigg(\tilde{L}^{\delta}+L_{f_1}^{\delta}t_{1}^{\delta}+L_{f_2}^{\delta}\frac{(t_{1})^{\alpha\delta}}{\alpha}\,\mathfrak{F}_{1}^{*}+L_{f_3}^{\delta}\frac{(t_1)^{\alpha\delta}}{\alpha}\,\mathfrak{F}_{2}^{*}\bigg):
i=1,\ldots,m\bigg\}<1. \label{eq0-th3}
\end{eqnarray}

\begin{proof} To obtain the desired result, we will consider the following application $\mathbf{F}:PC_{1-\gamma}(J,\mathbb{R})\rightarrow PC_{1-\gamma}(J,\mathbb{R})$ defined by
\begin{eqnarray*}
(\mathbf{F}\,u)(t) &=&\left\{ 
\begin{array}{l}
\xi _{i}(t,u(t)),\text{ }t\in (t_{i},s_{i}],\,\,\,\,\,\,\,\,i=1,2,\ldots ,m; \\ 
{\mathbb{\mathbb{P}}}_{\alpha ,\beta }(t-s)(u(0)-g(u))+\displaystyle\int_{0}^{t}\mathbf{K}%
_{\alpha }(t-s)\mathbf{F}_{\mathfrak{T},\mathfrak{V}}(s,u(s))ds,\,\,\,\,\,\,\,\,t\in
\lbrack 0,t_{1}]; \\ 
{\mathbb{\mathbb{P}}}_{\alpha ,\beta }(t-s)[\xi
_{i}(s_{i},u(s_{i}))-g(u)]+\displaystyle\int_{s_{i}}^{t}\mathbf{K}_{\alpha }(t-s)\mathbf{F}_{\mathfrak{T},\mathfrak{V}}(s,u(s))ds,\\ t\in
(s_{i},t_{i+1}],\,\,i=1.\dots ,m.%
\end{array}%
\right.  \\
(\mathbf{F}\,\mathcal{I}_{0^{+}}^{1-\gamma }u)(0) &=&u_{0}-g(u).
\end{eqnarray*}

Note that, $\mathbf{F}$ is well defined. Now, for any $u,v\in PC_{1-\gamma}(J,\mathbb{R})$ and $t\in(s_i,t_{i+1}]$, $i=1,2,\ldots,m$, we get
\begin{eqnarray}
\label{eq1-th3}
&&\bigg|t^{1-\gamma}[(\mathbf{F}\,u)(t)-(\mathbf{F}\,v)(t)]\bigg|\leq\bigg|t^{1-\gamma}{\mathbb{P}}_{\alpha,\beta}(t-s)
[\xi_{i}(s_i,u(s_i))-g(u)]\nonumber\\
&&+t^{1-\gamma}\int_{s_i}^{t}\mathbf{K}_{\alpha}(t-s)\mathbf{F}_{\mathfrak{T},\mathfrak{V}}(s,u(s))ds-t^{1-\gamma}{\mathbb{P}}_{\alpha,\beta}(t-s)[\xi_{i}(s_i,v(s_i))-g(v)]\nonumber\\
&&-t^{1-\gamma}\int_{s_i}^{t}\mathbf{K}_{\alpha}(t-s)\mathbf{F}_{\mathfrak{T},\mathfrak{V}}(s,v(s))ds\bigg|\nonumber\\
&&\leq\bigg|t^{1-\gamma}{\mathbb{P}}_{\alpha,\beta}(t-s)\bigg|\bigg|\xi_{i}(s_i,u(s_i))-\xi_{i}(s_i,v(s_i))\bigg|
+\bigg|t^{1-\gamma}{\mathbb{P}}_{\alpha,\beta}(t-s)\bigg|\bigg|g(u)-g(v)\bigg|\nonumber\\
&&+t^{1-\gamma}\int_{s_i}^{t}\mathbf{K}_{\alpha}(t-s)|\mathbf{F}_{\mathfrak{T},\mathfrak{V}}(s,u(s))-\mathbf{F}_{\mathfrak{T},\mathfrak{V}}(s,v(s))|ds\nonumber\\
&&\leq \mathfrak{M}\,L_{\xi_i}|u(s_i)-v(s_i)|+\mathfrak{M}\,\tilde{L}|u-v|+t^{1-\gamma}\int_{s_i}^{t}\mathfrak{M}L_{f_1}|u(s)-v(s)|ds \nonumber\\
&&+t^{1-\gamma}\int_{s_i}^{t}\mathfrak{M}L_{f_2}|\mathfrak{T}u(s)-\mathfrak{T}v(s)|ds+t^{1-\gamma}\int_{s_i}^{t}\mathfrak{M}L_{f_3}|\mathfrak{V}u(s)-\mathfrak{V}v(s)|ds.
\end{eqnarray}

Observe that,
\begin{eqnarray}
\label{eq2-th3}
&&t^{1-\gamma}\mathfrak{M}\int_{s_i}^{t}L_{f_2}|\mathfrak{T}u(s)-\mathfrak{T}v(s)|ds\nonumber\\
&&\leq t^{1-\gamma}\mathfrak{M}L_{f_2}\int_{s_i}^{t}\int_{0}^{s}|\mathbf{K}(s,\tau)|u(\tau)-v(\tau)|d\tau ds\nonumber\\
&&\leq K\int_{s_i}^{t}\sup_{t\in[s_i,t_{i+1}]}\bigg|t^{1-\gamma}[u(t)-v(t)]\bigg|(t-s)^{\alpha-1}
\sup_{t\in{J}}\int_{0}^{t}(s-\tau)^{1-\alpha}|\mathfrak{F}_{1}(s,\tau)|d\tau ds\nonumber\\
&&=\mathfrak{M}{L_{f_2}}\int_{s_i}^{t}{\lVert u-v \rVert}_{PC_{1-\gamma}}(t-s)^{\alpha-1}\mathfrak{F}_{1}^{*}ds\nonumber\\
&&\leq \mathfrak{M}{L_{f_2}}{\lVert u-v \rVert}_{PC_{1-\gamma}}\mathfrak{F}_{1}^{*}\frac{(t_{i+1}-s_i)^{\alpha}}{\alpha},
\end{eqnarray}
where $\gamma=\alpha+\beta(\alpha-1)$ with $0\leq \beta \leq 1$ and $0<\alpha\leq{1}$.

On the other hand, realizing the same steps of inequality \textnormal{Eq.(\ref{eq2-th3})}, we obtain 
\begin{eqnarray}
t^{1-\gamma}\mathfrak{M}\int_{s_i}^{b}L_{f_3}|\mathfrak{V}u(s)-\mathfrak{V}v(s)|ds\leq\mathfrak{M}\,L_{f_3}{\lVert u-v \rVert}_{PC_{1-\gamma}}
\mathfrak{F}_{2}^{*}\frac{(t_{i+1}-s_i)^{\alpha}}{\alpha}. \label{eq3-th3}
\end{eqnarray}

Substituting the inequality \textnormal{Eq.(\ref{eq2-th3})} and \textnormal{Eq.(\ref{eq3-th3})}, into inequality \textnormal{Eq.(\ref{eq1-th3})}, we get
\begin{eqnarray*}
&&\bigg|t^{1-\gamma}[(\mathbf{F}\,u)(t)-(\mathbf{F}\,v)(t)]\bigg|\\
&&\leq \mathfrak{M}\,L_{\xi_i}|u(s_i)-v(s_i)|+\mathfrak{M}\tilde{L}|u-v|+t^{1-\gamma}\int_{s_i}^{t}\mathfrak{M}
L_{f_1}|u(s)-v(s)|ds\\
&&+t^{1-\gamma}\int_{s_i}^{t}\mathfrak{M}L_{f_2}|\mathfrak{T}u(s)-\mathfrak{T}v(s)|ds+t^{1-\gamma}\int_{s_i}^{t}\mathfrak{M}
L_{f_3}|\mathfrak{V}u(s)-\mathfrak{V}v(s)|ds\\
&&\leq \mathfrak{M}\,{L}_{\xi_i}|u(s_i)-v(s_i)|_{PC_{1-\gamma}}+\mathfrak{M}\tilde{L}|u-v|_{PC_{1-\gamma}}
+\sup_{t\in{J}}\bigg|t^{1-\gamma}[u(t)-v(t)]\bigg|\int_{s_i}^{t}\mathfrak{M}\,L_{f_1}ds\\
&&+\mathfrak{M}{L_{f_2}}{\lVert u-v \rVert}_{PC_{1-\gamma}}\mathfrak{F}_{1}^{*}\frac{(t_{i+1}-s_i)^{\alpha}}{\alpha}
+\mathfrak{M}\,L_{f_3}{\lVert u-v \rVert}_{PC_{1-\gamma}}\mathfrak{F}_{2}^{*}\frac{(t_{i+1}-s_i)^{\alpha}}{\alpha}
\\&&\leq \mathfrak{M}\,L_{\xi_i}{\lVert u-v \rVert}_{PC_{1-\gamma}}+\mathfrak{M}\tilde{L}|u-v|_{PC_{1-\gamma}}
+\mathfrak{M}\,L_{f_1}(t_{i+1}-s_i){\lVert u-v \rVert}_{PC_{1-\gamma}}\\
&& +\mathfrak{M}\,L_{f_2}{\lVert u-v \rVert}_{PC_{1-\gamma}}\mathfrak{F}_{1}^{*}\frac{(t_{i+1}-s_i)^{\alpha}}{\alpha}
+\mathfrak{M}\,L_{f_3}{\lVert u-v \rVert}_{PC_{1-\gamma}}\mathfrak{F}_{2}^{*}\frac{(t_{i+1}-s_i)^{\alpha}}{\alpha}\\
&&\leq \mathfrak{M}(L_{\xi_i}^{\delta}+\tilde{L}^{\delta}){\lVert u-v \rVert}_{PC_{1-\gamma}}+
\mathfrak{M}\,L_{f_1}^{\delta}{\lVert u-v \rVert}_{PC_{1-\gamma}}(t_{i+1}-s_i)^{\delta}\\
&&+ \mathfrak{M}\bigg( L_{f_2}^{\delta}\mathfrak{F}_{1}^{*}\frac{(t_{i+1}-s_i)^{\alpha\delta}}{\alpha}
+L_{f_3}^{\delta}\mathfrak{F}_{2}^{*}\frac{(t_{i+1}-s_i)^{\alpha\delta}}{\alpha}\bigg) {\lVert u-v \rVert}_{PC_{1-\gamma}}
\end{eqnarray*}
which implies
\begin{eqnarray*}
&&\bigg|t^{1-\gamma}[(\mathbf{F}\,u)(t)-(\mathbf{F}\,v)(t)]\bigg|^{\delta}\leq 
\mathfrak{M}(L_{\xi_i}^{\delta}+\tilde{L}^{\delta}){\lVert u-v \rVert}_{PC_{1-\gamma}}+
\mathfrak{M}\,L_{f_1}^{\delta}{\lVert u-v \rVert}_{PC_{1-\gamma}}(t_{i+1}-s_i)^{\delta}\\
&&+ \mathfrak{M}\bigg( L_{f_2}^{\delta}\mathfrak{F}_{1}^{*}\frac{(t_{i+1}-s_i)^{\alpha\delta}}{\alpha}
+L_{f_3}^{\delta}\mathfrak{F}_{2}^{*}\frac{(t_{i+1}-s_i)^{\alpha\delta}}{\alpha}\bigg) {\lVert u-v \rVert}_{PC_{1-\gamma},\delta}.
\end{eqnarray*}

In this sense, we obtain
\begin{eqnarray}
{\lVert \mathbf{F}\,u-\mathbf{F}\,v \rVert}_{\mathcal{PC}_{1-\gamma,\delta}}&\leq& \mathfrak{M}(L_{\xi_i}^{\delta}+\tilde{L}^{\delta})
{\lVert u-v \rVert}_{PC_{1-\gamma},\delta}+\mathfrak{M}\bigg(L_{f_1}^{\delta}(t_{i+1}-s_i)^{\delta}+
L_{f_2}^{\delta}\frac{(t_{i+1}-s_i)^{\alpha\delta}}{\alpha}\mathfrak{F}_{1}^{*}\nonumber\\
&+&L_{f_3}^{\delta}\frac{(t_{i+1}-s_i)^{\alpha\delta}}{\alpha}\mathfrak{F}_{2}^{*}\bigg){\lVert u-v \rVert}_{PC_{1-\gamma},\delta},
\quad t\in(s_i,t_{i+1}]. \label{eq4-th3}
\end{eqnarray}

On the other hand, realizing the procedure as in the inequality \textnormal{Eq.(\ref{eq4-th3})}, we have
\begin{eqnarray}
&& {\lVert \mathbf{F}\,u-\mathbf{F}\,v \rVert}_{\mathcal{PC}_{1-\gamma,\delta}}\leq \mathfrak{M}\tilde{L}^{\delta}
{\lVert u-v \rVert}_{PC_{1-\gamma},\delta}+\mathfrak{M}\bigg(L_{f_1}^{\delta}(t_{1})^{\delta}+
L_{f_2}^{\delta}\frac{(t_{1})^{\alpha\delta}}{\alpha}\mathfrak{F}_{1}^{*}\nonumber\\
&&+L_{f_3}^{\delta}\frac{(t_{1})^{\alpha\delta}}{\alpha}\mathfrak{F}_{2}^{*}\bigg){\lVert u-v \rVert}_{PC_{1-\gamma},\delta},
\quad t\in[0,1] \label{eq5-th3}
\end{eqnarray}
and
\begin{eqnarray}
{\lVert \mathbf{F}\,u-\mathbf{F}\,v \rVert}_{\mathcal{PC}_{1-\gamma,\delta}}\leq(\mathfrak{M}\tilde{L}^{\delta}+
\mathfrak{M}L_{\xi_i}^{\delta}){\lVert u-v \rVert}_{PC_{1-\gamma},\delta}, \quad t\in(t_i,s_i].\label{eq6-th3}
\end{eqnarray}

From the inequalities \textnormal{Eq.(\ref{eq4-th3})}, \textnormal{Eq.(\ref{eq5-th3})} and
\textnormal{Eq.(\ref{eq6-th3})}, we have
\begin{eqnarray*}
{\lVert \mathbf{F}\,u-\mathbf{F}\,v \rVert}_{\mathcal{PC}_{1-\gamma,\delta}}\leq\Lambda\,{\lVert u-v \rVert}_{PC_{1-\gamma},\delta},
\end{eqnarray*}
where
\begin{eqnarray*}
\Lambda&:=&max\bigg\{\mathfrak{M}({L}_{\xi_i}^{\delta}+\tilde{L}^{\delta})+\mathfrak{M}
\bigg(L_{f_1}^{\delta}(t_{i+1}-s_i)^{\delta}+L_{f_2}^{\delta}\frac{(t_{i+1}-s_i)^{\alpha\delta}}{\alpha}\,\mathfrak{F}_{1}^{*}
+L_{f_3}^{\delta}\frac{(t_{i+1}-s_i)^{\alpha\delta}}{\alpha}\,\mathfrak{F}_{2}^{*}\bigg);\\
&&\mathfrak{M}\bigg(\tilde{L}^{\delta}+L_{f_1}^{\delta}t_{1}^{\delta}+L_{f_2}^{\delta}\frac{(t_{1})^{\alpha\delta}}{\alpha}\,\mathfrak{F}_{1}^{*}+L_{f_3}^{\delta}\frac{(t_1)^{\alpha\delta}}{\alpha}\,\mathfrak{F}_{2}^{*}\bigg):
i=1,\ldots,m\bigg\}<1.
\end{eqnarray*}

Thus, we conclude that $\mathbf{F}$ is a contraction mapping. As $PC_{1-\gamma}(J,\mathbb{R})$ is a Banach space endowed with the norm ${\lVert\, \cdot \,\rVert}_{PC_{1-\gamma},\delta}$ and $\mathbf{F}: {PC_{1-\gamma}}(J,\mathbb{R})\rightarrow{PC_{1-\gamma}}(J,\mathbb{R})$ is a contraction mapping, then by Banach fixed point theorem, the \textnormal{Eq.(\ref{eq3-th3})} has a unique solution.
\end{proof}
\end{teorema}

Next we will investigate the stability $\delta$-Ulam-Hyers-Rassias of the solution of Eq.(\ref{eq:1}).

\begin{teorema} Assume that \textnormal{(H1)-(H5)} and \textnormal{Eq.(\ref{eq0-th3})} are satisfied. Then \textnormal{Eq.(\ref{eq:1})} is generalized $\delta$-Ulam-Hyers-Rassias stable with respect to $(\varphi,\phi)$.

\begin{proof}
Let $v\in PC_{1-\gamma}^{1}(J,\mathbb{R})$ be a solution of \textnormal{Eq.(\ref{ineq-2})} and $u$ the unique solution of the impulsive fractional Cauchy problem, given by
\begin{eqnarray*}
\left\{
\begin{array}{lllllll}
^{\mathbf{H}}\mathfrak{D}_{0^+}^{\alpha,\beta}u(t)=\mathbf{F}_{\mathfrak{T},\mathfrak{V}}(t,u(t)),\,\,\,\,\,\,\,\, t\in(s_i,t_{i+1}],\,\, i\in[0,m] \\
u(t)=\xi_i(t,u(t)), t\in(t_i,s_i], \,\,\,\,\,\,\,\,i\in[1,m]\\
I_{0^+}^{1-\gamma}u(0)+g(u)=u_0. \\
\end{array}
\right.
\end{eqnarray*}

In this way, we have
\begin{eqnarray*}
u(t)=
\left\{
\begin{array}{lllllll}
{\mathbb{P}}_{\alpha,\beta}(t-s)[\xi_i(t,u(t))-g(t)],\,\, t\in(t_i,s_{i}],\,\, i\in[1,m] \\
\displaystyle {\mathbb{P}}_{\alpha,\beta}(t-s)[u(0)-g(u)]+\int_{0}^{t}\mathbf{K}_{\alpha}(t-s)\mathbf{F}_{\mathfrak{T},\mathfrak{V}}(s,u(s)), t\in[0,t_1],\\
\displaystyle {\mathbb{P}}_{\alpha,\beta}(t-s)[\xi_i(s_i,u(s_i))-g(u)]+\int_{s_i}^{t}\mathbf{K}_{\alpha}(t-s)\mathbf{F}_{\mathfrak{T},\mathfrak{V}}(s,u(s))ds,\\
t\in[s_i,t_{i+1}],\quad i=[1,m].\\
\end{array}
\right.
\end{eqnarray*}

Since \textnormal{Eq.(\ref{ineq-rk})} holds, for each $t\in(s_i,t_{i+1}]$, $i=1,2,\ldots,m$, we get
\begin{eqnarray*}
&&\bigg|v(t)-{\mathbb{P}}_{\alpha,\beta}(t-s)[\xi_i(s_i,v(s_i))-g(v)]-\int_{s_i}^{t}\mathbf{K}_{\alpha}(t-s)\mathbf{F}_{\mathfrak{T},\mathfrak{V}}(s,v(s))ds\bigg|\\
&&\leq {\lVert {\mathbb{P}}_{\alpha,\beta}(t-s) \rVert}_{PC_{1-\gamma}}\phi+\int_{s_i}^{t}{\lVert \mathbf{K}_{\alpha}(t-s) \rVert}_{PC_{1-\gamma}}\varphi(s)ds\\
&&\leq\mathfrak{M}\phi+\mathfrak{M}\int_{s_i}^{t}\varphi(s)ds\leq\mathfrak{M}\,\phi+\mathfrak{M}c_{\varphi}\varphi(t),
\end{eqnarray*}
and for $t\in(t_i,s_i]$, $i=1,2,\ldots,m$, we get
\begin{eqnarray*}
\bigg|v(t)-\mathbb{P}_{\alpha,\beta}(t-s)[\xi_i(s_i,v(s_i))-g(v)]\bigg|\leq\mathfrak{M}\,\phi
\end{eqnarray*}
and for $t\in[0,t_1]$,
\begin{eqnarray*}
&&\bigg|v(t)-{\mathbb{P}}_{\alpha,\beta}(t-s)[v(0)-g(v)]\int_{0}^{t}\mathbf{K}_{\alpha}(t-s)\mathbf{F}_{\mathfrak{T},\mathfrak{V}}(s,v(s))ds\bigg|\\
&&\leq\int_{0}^{t} {\lVert \mathbf{K}_{\alpha}(t-s) \rVert}_{PC_{1-\gamma}}\varphi(s)ds
\leq\mathfrak{M}\int_{0}^{t}\varphi(s)ds\leq\mathfrak{M}c_{\varphi}\varphi(t).
\end{eqnarray*}

Thus, for each $t\in(s_i,t_{i+1}]$, $i=1,2,\ldots,m$, we obtain
\begin{eqnarray*}
&&|v(t)-u(t)|=\bigg| v(t)-\mathfrak{M}\,\xi_i(s_i,u(s_i))+\mathfrak{M}g(u)+\mathfrak{M}\int_{s_i}^{t}
\mathbf{F}_{\mathfrak{T},\mathfrak{V}}(s,v(s))ds\bigg|\\
&&\leq\bigg| v(t)-\mathfrak{M}\,\xi_i(s_i,v(s_i))+\mathfrak{M}g(v)-\mathfrak{M}\int_{s_i}^{t}\mathbf{F}_{\mathfrak{T},\mathfrak{V}}(s,v(s))ds\bigg|\\
&&+\mathfrak{M}|g(v)-g(u)|+\mathfrak{M}|\xi_i(s_i,v(s_i))-\xi_i(s_i,u(s_i))|\\
&&+ \mathfrak{M}\bigg(\int_{s_i}^{t}[\mathbf{F}_{\mathfrak{T},\mathfrak{V}}(s,v(s))-\mathbf{F}_{\mathfrak{T},\mathfrak{V}}(s,u(s))]ds\bigg)\\
&&\leq\mathfrak{M}(1+c_{\varphi})(\phi+\varphi(t))+\mathfrak{M}\tilde{L}|v(s)-u(s)|+\mathfrak{M}\sum_{0<s_i<t}
L_{\xi_i}|v(s_i)-u(s_i)|\\
&+& \mathfrak{M}\bigg(\int_{s_i}^{t}[L_{f_1}|v(s)-u(s)|+L_{f_2}\sup_{s\in{J}}|\mathfrak{T}v(s)-\mathfrak{T}u(s)|+L_{f_3}\sup_{s\in{J}}
|\mathfrak{V}v(s)-\mathfrak{V}u(s)|]ds\bigg)\\
&&\leq{\mathfrak{M}}(1+c_{\varphi})(\phi+\varphi(t))+\mathfrak{M}\tilde{L}|v(s)-u(s)|+\mathfrak{M}\sum_{0<s_i<t}
L_{\xi_i}|v(s_i)-u(s_i)|\\
&&+\mathfrak{M}\int_{0}^{t}L_{f_i}(\mathfrak{F}_{3}^{*}+\mathfrak{F}_{1}^{*}+\mathfrak{F}_{2}^{*})(t-s)^{\alpha-1}
|v(s)-u(s)|ds.
\end{eqnarray*}

Now, choosing $\tilde{\Phi}(t):=\mathfrak{M}(1+c_{\varphi})(\phi+\varphi(t))$, $\tilde{\delta}=\mathfrak{M}\tilde{L}$, $\beta_{k}=\mathfrak{M}L_{\xi_i}$, $g(t)=\mathfrak{M}L_{f_i} (\mathfrak{F}_{3}^{*}+\mathfrak{F}_{1}^{*}+\mathfrak{F}_{2}^{*})$ and
$t\in(s_i,t_{i+1}]$ and it is clear that both functions are nondecreasing and $\tilde{\Phi}, \tilde{\delta},g\in PC_{1-\gamma}(I,\mathbb{R})$.

Then, by \textnormal{Corollary \ref{Cor-A}}, we have
\begin{eqnarray*}
&&|v(t)-u(t)|\\
&\leq&\tilde{\Phi}(t)\bigg[\tilde{\delta}\,\mathbb{E}[g(t)\Gamma(\alpha)t^{\alpha}]\bigg(1+\beta\,\mathbb{E}_{\alpha}[g(t)\Gamma(\alpha)t^{\alpha}]\bigg)^{k-1}+\bigg(1+\beta\,\mathbb{E}_{\alpha}[g(t)\Gamma(\alpha)t^{\alpha}]\bigg)^{k}\bigg]
\mathbb{E}_{\alpha}[g(t)\Gamma(\alpha)t^{\alpha}]\\
&=& \tilde{\Phi}(t)\bigg[\mathfrak{M}\tilde{L}\,\mathbb{E}_{\alpha}[\mathfrak{M}L_{f_i}(\mathfrak{F}_{3}^{*}+\mathfrak{F}_{1}^{*}+\mathfrak{F}_{2}^{*})\Gamma(\alpha)t^{\alpha}]\bigg(1+\mathfrak{M}L_{\xi}\,\mathbb{E}_{\alpha}[\mathfrak{M}L_{f_i}(\mathfrak{F}_{3}^{*}+\mathfrak{F}_{1}^{*}+\mathfrak{F}_{2}^{*})\Gamma(\alpha)t^{\alpha}]\bigg)^{k}\bigg]\\
&&\times \mathbb{E}_{\alpha}[\mathfrak{M}L_{f_i}(\mathfrak{F}_{3}^{*}+\mathfrak{F}_{1}^{*}+\mathfrak{F}_{2}^{*})\Gamma(\alpha)t^{\alpha}],
\end{eqnarray*}
where $L_{\xi}=\max\{L_{\xi_1},L_{\xi_2},\ldots,L_{\xi_m}\}$. Thus,
\begin{eqnarray}
&&|v(t)-u(t)|^{\delta}\leq\bigg\{\tilde{\Phi}(t)\bigg[\mathfrak{M}\tilde{L}\,\mathbb{E}_{\alpha}[\mathfrak{M}L_{f_i}(\mathfrak{F}_{3}^{*}+\mathfrak{F}_{1}^{*}+\mathfrak{F}_{2}^{*})\Gamma(\alpha)t^{\alpha}]+\nonumber\\
&&\bigg(1+\mathfrak{M}L_{\xi}\,\mathbb{E}_{\alpha}[\mathfrak{M}L_{f_i}(\mathfrak{F}_{3}^{*}+\mathfrak{F}_{1}^{*}+\mathfrak{F}_{2}^{*})\Gamma(\alpha)t^{\alpha}]\bigg)^{k-1}\nonumber\\
&&+\bigg(1+\mathfrak{M}L_{\xi}\,\mathbb{E}_{\alpha}[\mathfrak{M}L_{f_i}(\mathfrak{F}_{3}^{*}+\mathfrak{F}_{1}^{*}+\mathfrak{F}_{2}^{*})\Gamma(\alpha)t^{\alpha}]\bigg)^{k}\,\bigg]\mathbb{E}_{\alpha}[\mathfrak{M}L_{f_i}(\mathfrak{F}_{3}^{*}+\mathfrak{F}_{1}^{*}+\mathfrak{F}_{2}^{*})\Gamma(\alpha)t^{\alpha}]\bigg\}^{\delta}\nonumber\\
&&\leq\bigg\{\mathfrak{M}(\phi+c_{\varphi})\bigg[\mathfrak{M}\tilde{L}\,\mathbb{E}_{\alpha}[\mathfrak{M}L_{f_i}(\mathfrak{F}_{3}^{*}+\mathfrak{F}_{1}^{*}+\mathfrak{F}_{2}^{*})\Gamma(\alpha){t^{\alpha}_{i+1}}]\nonumber\\
&&+\bigg(1+\mathfrak{M}L_{\xi}\,
\mathbb{E}_{\alpha}[\mathfrak{M}L_{f_i}(\mathfrak{F}_{3}^{*}+\mathfrak{F}_{1}^{*}+\mathfrak{F}_{2}^{*})\Gamma(\alpha){t^{\alpha}_{i+1}}]\bigg)^{k-1}+\mathbb{E}_{\alpha}[\mathfrak{M}L_{f_i}(\mathfrak{F}_{3}^{*}+\mathfrak{F}_{1}^{*}+\mathfrak{F}_{2}^{*})\Gamma(\alpha){t^{\alpha}_{i+1}}]\bigg)^{k}\,\bigg]\nonumber\\
&&\times \mathbb{E}_{\alpha}[\mathfrak{M}L_{f_i}(\mathfrak{F}_{3}^{*}+\mathfrak{F}_{1}^{*}+\mathfrak{F}_{2}^{*})\Gamma(\alpha){t^{\alpha}_{i+1}}]\bigg\}^{\delta}(\phi^{\delta}+\varphi(t)^{\delta}), \quad t\in(s_i,t_{i+1}],\quad i=1,2,\ldots,m.
\label{eq1-th4}
\end{eqnarray}

In addition, for $t\in(t_i,s_i]$, $i=1,2,\ldots,m$, we get
\begin{eqnarray}
|v(t)-u(t)|&=&\bigg|v(t)-{\mathbb{P}}_{\alpha,\beta}(t-s)\bigg(\xi_{i}(t,u(t))-g(t)\bigg)\bigg|\nonumber\\
&=&\bigg|v(t)-{\mathbb{P}}_{\alpha,\beta}(t-s)\xi_{i}(t,u(t))+{\mathbb{P}}_{\alpha,\beta}(t-s)g(t)\bigg|\nonumber\\
&\leq& \mathfrak{M}|v(t)-\xi_{i}(t,u(t))+\mathfrak{M}g(u)|^{\delta}.
\label{eq2-th4}
\end{eqnarray}

Taking $\delta$ on both sides of the inequality \textnormal{Eq.(\ref{eq2-th4})}, we have
\begin{eqnarray*}
|v(t)-u(t)|^{\delta}&\leq&\mathfrak{M}^{\delta}|v(t)-\xi_{i}(t,u(t))+\mathfrak{M}g(u)|^{\delta}
\leq\mathfrak{M}|v(t)-\xi_{i}(t,u(t))+\mathfrak{M}g(u)|^{\delta}\nonumber\\
&\leq&\mathfrak{M}|v(t)-\xi_{i}(t,v(t))+g(v)|^{\delta}+\mathfrak{M}|g(v)-g(u)|^{\delta}+\mathfrak{M}|\xi_{i}(t,v(t))-
\xi_{i}(t,u(t))|^{\delta}\nonumber\\
&\leq&\mathfrak{M}\phi^{\delta}+\mathfrak{M}\tilde{L}^{\delta}|v(t)-u(t)|^{\delta}+\mathfrak{M}L_{\xi_i}^{\delta}|v(t)-u(t)|^{\delta}\nonumber\\
&&\Rightarrow (1-\mathfrak{M}\tilde{L}^{\delta}-\mathfrak{M}L_{\xi_i}^{\delta})|v(t)-u(t)|^{\delta}\leq\mathfrak{M}\phi^{\delta}.
\end{eqnarray*}

This implies,
\begin{eqnarray}
|v(t)-u(t)|^{\delta}\leq \frac{\mathfrak{M}\phi^{\delta}}{1-\mathfrak{M}\tilde{L}^{\delta}-\mathfrak{M}L_{\xi_i}^{\delta}}.\label{eq3-th4}
\end{eqnarray}

On the other hand, for $t\in[0,t_1]$, we have
\begin{eqnarray*}
&&|v(t)-u(t)|=\bigg|v(t)-\mathfrak{M}(u(0)-g(u))-\mathfrak{M}\int_{0}^{t}\mathbf{F}_{\mathfrak{T},\mathfrak{V}}(s,u(s))ds\bigg|\nonumber\\
&=&\bigg|v(t)-\mathfrak{M}(v(0)-g(v))-\mathfrak{M}\int_{0}^{t}\mathbf{F}_{\mathfrak{T},\mathfrak{V}}(s,v(s))ds+\nonumber\\
&+&\mathfrak{M}(v(0)-g(v))-\mathfrak{M}(u(0)-g(u))-\mathfrak{M}\int_{0}^{t}[\mathbf{F}_{\mathfrak{T},\mathfrak{V}}(s,u(s))-\mathbf{F}_{\mathfrak{T},\mathfrak{V}}(s,v(s))]ds\bigg|
\nonumber\\
&\leq&\mathfrak{M}c_{\varphi}\varphi(t)+\mathfrak{M}\tilde{L}|v(t)-u(t)|+\mathfrak{M}|v(0)-u(0)|
+\mathfrak{M}\int_{0}^{t}L_{\xi_i}(\mathfrak{F}_{3}^{*}+\mathfrak{F}_{1}^{*}+\mathfrak{F}_{2}^{*})(t-s)^{\alpha-1}
|v(s)-u(s)|ds\nonumber\\
&\leq&\mathfrak{M}c_{\varphi}\varphi(t)+(\mathfrak{M}\tilde{L}+\mathfrak{M})|v(t)-u(t)|
+\mathfrak{M}\int_{0}^{t}L_{\xi_i}(\mathfrak{F}_{3}^{*}+\mathfrak{F}_{1}^{*}+\mathfrak{F}_{2}^{*})(t-s)^{\alpha-1}
|v(s)-u(s)|ds.
\end{eqnarray*}

By \textnormal{Theorem \ref{theo-1}}, we have
\begin{eqnarray*}
&&|v(t)-u(t)|\\
&\leq&\mathfrak{M}c_{\varphi}\varphi(t)\bigg[(\mathfrak{M}\tilde{L}+\mathfrak{M})\mathbb{E}_{\alpha}
[\mathfrak{M}L_{\xi_i}(\mathfrak{F}_{3}^{*}+\mathfrak{F}_{1}^{*}+\mathfrak{F}_{2}^{*})\Gamma(\alpha)t^{\alpha}]+1\bigg]
\mathbb{E}_{\alpha}[\mathfrak{M}L_{\xi_i}(\mathfrak{F}_{3}^{*}+\mathfrak{F}_{1}^{*}+\mathfrak{F}_{2}^{*})\Gamma(\alpha)t^{\alpha}]\nonumber\\
&=&\mathfrak{M}(\mathfrak{M}+\tilde{L})c_{\varphi}\varphi(t)\bigg[(\mathfrak{M}\tilde{L}+\mathfrak{M})\mathbb{E}_{\alpha}
[\mathfrak{M}L_{\xi_i}(\mathfrak{F}_{3}^{*}+\mathfrak{F}_{1}^{*}+\mathfrak{F}_{2}^{*})\Gamma(\alpha)t_{1}^{\alpha}]+1\bigg]
\mathbb{E}_{\alpha}[\mathfrak{M}L_{\xi_i}(\mathfrak{F}_{3}^{*}+\mathfrak{F}_{1}^{*}+\mathfrak{F}_{2}^{*})\Gamma(\alpha)t_{1}^{\alpha}].
\end{eqnarray*}

In this sense, we obtain
\begin{eqnarray}
&&|v(t)-u(t)|^{\delta}\nonumber\\
&\leq&\bigg\{\mathfrak{M}(\mathfrak{M}+\tilde{L})c_{\varphi}\varphi(t)\bigg[(\mathfrak{M}\tilde{L}+\mathfrak{M})\mathbb{E}_{\alpha}[\mathfrak{M}L_{\xi_i}(\mathfrak{F}_{3}^{*}+\mathfrak{F}_{1}^{*}+\mathfrak{F}_{2}^{*})\Gamma(\alpha)t_{1}^{\alpha}]+1\bigg]
\nonumber\\
&\times& \mathbb{E}_{\alpha}[\mathfrak{M}L_{\xi_i}(\mathfrak{F}_{3}^{*}+\mathfrak{F}_{1}^{*}+\mathfrak{F}_{2}^{*})\Gamma(\alpha)t_{1}^{\alpha}]\bigg\}^{\delta}\nonumber\\
&\leq&\bigg\{\mathfrak{M}(\mathfrak{M}+\tilde{L})c_{\varphi}\varphi(t)\bigg[(\mathfrak{M}\tilde{L}+\mathfrak{M})\mathbb{E}_{\alpha}[\mathfrak{M}L_{\xi_i}(\mathfrak{F}_{3}^{*}+\mathfrak{F}_{1}^{*}+\mathfrak{F}_{2}^{*})\Gamma(\alpha)t_{1}^{\alpha}]+1\bigg]\nonumber\\
&\times& \mathbb{E}_{\alpha}[\mathfrak{M}L_{\xi_i}(\mathfrak{F}_{3}^{*}+\mathfrak{F}_{1}^{*}+\mathfrak{F}_{2}^{*})\Gamma(\alpha)t_{1}^{\alpha}]\bigg\}^{\delta}\varphi^{\delta}(t), \quad t\in[0,t_1]. \label{eq4-th4}
\end{eqnarray}

By means of the \textnormal{Eq.(\ref{eq1-th4})}, \textnormal{Eq.(\ref{eq3-th4})} and  \textnormal{Eq.(\ref{eq4-th4})}, we get
\begin{eqnarray*}
&&|v(t)-u(t)|^{\delta}\\
&\leq&\bigg(\bigg\{\mathfrak{M}(1+c_{\varphi})\bigg[\mathfrak{M}\tilde{L}
\mathbb{E}_{\alpha}[\mathfrak{M}L_{f_i}(\mathfrak{F}_{3}^{*}+\mathfrak{F}_{1}^{*}+\mathfrak{F}_{2}^{*})\Gamma(\alpha)t_{i+1}^{\alpha}]\\
&&\times \bigg(1+\mathfrak{M}L_{\xi}\,\mathbb{E}_{\alpha}[\mathfrak{M}L_{f_i}(\mathfrak{F}_{3}^{*}+\mathfrak{F}_{1}^{*}+\mathfrak{F}_{2}^{*})\Gamma(\alpha)t_{i+1}^{\alpha}]\bigg)^{k-1}+\bigg(1+\mathfrak{M}L_{\xi}\,\mathbb{E}_{\alpha}[\mathfrak{M}L_{f_i}(\mathfrak{F}_{3}^{*}+\mathfrak{F}_{1}^{*}+\mathfrak{F}_{2}^{*})\Gamma(\alpha)t_{i+1}^{\alpha}]\bigg)^{k}\,\bigg]\\
&&\times\mathbb{E}_{\alpha}[\mathfrak{M}L_{f_i}(\mathfrak{F}_{3}^{*}+\mathfrak{F}_{1}^{*}+\mathfrak{F}_{2}^{*})\Gamma(\alpha)t_{i+1}^{\alpha}]
\bigg\}^{\delta}+\frac{\mathfrak{M}\phi^{\delta}}{1-\mathfrak{M}\tilde{L}^{\delta}-\mathfrak{M}L_{f_i}^{\delta}}\\
&&+\mathfrak{M}(\mathfrak{M}+\tilde{L})c_{\varphi}\bigg[\mathbb{E}_{\alpha}[\mathfrak{M}L_{\xi_i}(\mathfrak{F}_{3}^{*}+\mathfrak{F}_{1}^{*}+\mathfrak{F}_{2}^{*})\Gamma(\alpha)t_{1}^{\alpha}]+1\bigg]\mathbb{E}_{\alpha}[\mathfrak{M}L_{f_i}(\mathfrak{F}_{3}^{*}+\mathfrak{F}_{1}^{*}+\mathfrak{F}_{2}^{*})\Gamma(\alpha)t_{1}^{\alpha}](\phi^{\delta}+\varphi^{\delta}(t))\\
&&:=C_{f_i,\xi_i,\varphi}^{\mathfrak{M},\tilde{L},\delta}(\phi^{\delta}+\varphi^{\delta}(t)), \quad t\in{J},
\end{eqnarray*}
which implies that \textnormal{Eq.(\ref{eq:1})} is generalized $\delta$-Ulam-Hyers-Rassias stable with respect to $(\varphi, \phi)$.
\end{proof}
\end{teorema}

\section{Concluding Remarks}

We conclude this paper with the objectives achieved, i.e., we present a generalization for the Gronwall inequality generating a new class of inequalities of the Gronwall type as well as other generalizations, allowing new tools in the study of Ulam-Hyers stability of fractional differential equations. On the other hand, the second major result of the paper, the investigation of the uniqueness and $\delta$-Ulam-Hyers-Rassias stability of the impulsive differential equation with non-instantaneous impulses by means of the Hilfer fractional derivative, was obtained.

An interesting observation to be point out, we use Hilfer fractional derivative to perform this paper. However, the following question emerges: Why we don't use the $\psi$-Hilfer fractional derivative, since it contains as particular cases, a wide class of fractional derivatives, in particular, the Hilfer fractional derivative? The main reason for not using the $\psi$-Hilfer fractional derivative is the fact that the Laplace transform is not available with respect to other functions to obtain mild solutions for fractional differential equations formulated by means of the $\psi$-Hilfer fractional derivative. This is an open problem in the field of fractional calculus. The second case that prevents us from carrying out the studies here, is a Leibniz rule for the $\psi$-Hilfer fractional derivative. However, recently Sousa and Oliveira, proposed the so-called Leibniz type rule \cite{rule}.

To conclude, we leave here an open door to the study of existence, uniqueness and Ulam-Hyers stabilities for the various types of differential equations formulated by means of the $\psi$-Hilfer fractional derivative.

\section{Acknowledgements}
I  have  been  financially  supported  by  PNPD-CAPES  scholarship  of  the  Pos-Graduate Program in Applied Mathematics IMECC-Unicamp.

\bibliography{ref}
\bibliographystyle{plain}

\end{document}